\documentclass[12pt]{article}

\usepackage{amsmath,amsthm,amsfonts,amssymb,latexsym,amscd}

\font\fiverm=cmr5

\input{prepictex.tex}
\input{pictex.tex}
\input{postpictex.tex}

\begin{document}

\newtheorem*{theo}{Theorem}
\newtheorem*{pro}{Proposition}
\newtheorem*{cor}{Corollary}
\newtheorem*{lem}{Lemma}
\newtheorem{theorem}{Theorem}[section]
\newtheorem{corollary}[theorem]{Corollary}
\newtheorem{lemma}[theorem]{Lemma}
\newtheorem{proposition}[theorem]{Proposition}
\newtheorem{conjecture}[theorem]{Conjecture}
\newtheorem{definition}[theorem]{Definition}
\newtheorem{problem}[theorem]{Problem}
\newtheorem{remark}[theorem]{Remark}
\newtheorem{example}[theorem]{Example}
\newcommand{\Naturali}{{\mathbb{N}}}
\newcommand{\Reali}{{\mathbb{R}}}
\newcommand{\Complessi}{{\mathbb{C}}}
\newcommand{\Toro}{{\mathbb{T}}}
\newcommand{\Relativi}{{\mathbb{Z}}}
\newcommand{\HH}{\mathfrak H}
\newcommand{\KK}{\mathfrak K}
\newcommand{\LL}{\mathfrak L}
\newcommand{\as}{\ast_{\sigma}}
\newcommand{\tn}{\vert\hspace{-.3mm}\vert\hspace{-.3mm}\vert}
\def\mA{{\mathfrak A}}
\def\A{{\mathcal A}}
\def\mB{{\mathfrak B}}
\def\B{{\mathcal B}}
\def\C{{\mathcal C}}
\def\D{{\mathcal D}}
\def\F{{\mathcal F}}
\def\H{{\mathcal H}}
\def\J{{\mathcal J}}
\def\K{{\mathcal K}}
\def\L{{\cal L}}
\def\N{{\cal N}}
\def\M{{\cal M}}
\def\O{{\mathcal O}}
\def\P{{\cal P}}
\def\S{{\cal S}}
\def\T{{\cal T}}
\def\U{{\cal U}}
\def\W{{\cal W}}
\def\b{\lambda_B(P}
\def\j{\lambda_J(P}

\title{Labeled Trees and Localized Automorphisms of the Cuntz Algebras}

\author{Roberto Conti, Wojciech Szyma{\'n}ski}

\date{}
\maketitle

\renewcommand{\sectionmark}[1]{}

\noindent
{\small \date{10 May, 2008}}

\vspace{7mm}
\begin{abstract}
We initiate a detailed and systematic study of automorphisms of the Cuntz algebras $\O_n$
which preserve both the diagonal and the core $UHF$-subalgebra. A general criterion
of invertibility of endomorphisms yielding such automorphisms is given.
Combinatorial investigations of endomorphisms related to permutation matrices are presented.
Key objects entering this analysis are labeled rooted trees equipped with additional data.
Our analysis provides insight into the structure of ${\rm Aut}(\O_n)$
and leads to numerous new examples. In particular, we completely classify all
such automorphisms of ${\mathcal O}_2$ for the permutation unitaries in $\otimes^4 M_2$.
We show that the subgroup of ${\rm Out}(\O_2)$ generated by these automorphisms contains
a copy of the infinite dihedral group ${\mathbb Z} \rtimes {\mathbb Z}_2$.
\end{abstract}

\vfill
\noindent {\bf MSC 2000}: 46L40, 46L05, 37B10

\vspace{3mm}
\noindent {\bf Keywords}: Cuntz algebra, endomorphism, automorphism,
  Cartan subalgebra, core $UHF$-subalgebra, normalizer, permutation, tree.

\newpage

\vspace*{40mm}

\noindent
`Nel mezzo del cammin di nostra vita\\
mi ritrovai per una selva oscura,\\
ch\'e la diritta via era smarrita.'

\vspace{3mm}\noindent
Dante Alighieri, {\em La Divina Commedia, Inferno}

\vspace{16mm}\noindent
`Macbeth shall never vanquished be, until\\
great Birnam wood to high Dunsinane hill\\
shall come against him.'

\vspace{3mm}\noindent
William Shakespeare, {\em Macbeth}

\newpage

\tableofcontents

\section{Introduction}

In recent years endomorphisms of Cuntz algebras have received a lot of
attention and have been deeply investigated from the point of view of
Jones index theory and sector theory \cite{i,j,l,CP,JS,CF,Go,Kaw1,Kaw2}.
In these theories, emphasis is often placed on proper endomorphisms rather
than automorphisms. However, automorphisms
of Cuntz algebras have also been studied, sometimes in connection with
classification of group actions (for example, see \cite{Arc,Cun2,MaTo,KK,i2,Matui}).
In the present paper, our main interest lies in classification of a special
class of localized automorphisms of ${\mathcal O}_n$.

In his beautiful paper \cite{Cun2}, Joachim Cuntz
initiated systematic investigations of the automorphism group of $\O_n$. In particular, he
showed that the group of those automorphisms which preserve the diagonal subalgebra
contains a maximal abelian normal subgroup whose quotient (the Weyl group) is discrete.
Restricting even further to those automorphisms which preserve both the diagonal and the
core $UHF$-subalgebra one finds even nicer structure. Thus Cuntz suggested
that classification of all elements of this restricted Weyl group `is a combinatorial
problem, and should be possible'. By now thirty years have passed and this
classification has not been achieved, nor is it even in sight. Presumably,
this is due not to lack of interest of high power researchers in this exciting problem
but rather to great difficulties involved. It appears that there are two sources
of these difficulties.

Firstly, as demonstrated by Cuntz \cite{Cun2}, automorphisms of $\O_n$
are best seen as special class of endomorphisms. The latter are in
a one-to-one correspondence with unitary elements of $\O_n$ via a certain natural
correspondence ${\mathcal U}(\O_n)\ni u \mapsto \lambda_u \in {\rm End}(\O_n)$.
The problem is that in general there is no easy way of verifying which
unitaries $u$ give rise to invertible endomorphisms (i.e. automorphisms) $\lambda_u$.
In the present article we provide a remedy to this problem for a large class of
endomorphisms related to unitary matrices in $M_{n^k}({\mathbb C})$ contained
in the $UHF$-subalgebra (see Theorem \ref{main}, below).

Secondly, again as shown by Cuntz \cite{Cun2}, analysis of the Weyl group reduces
to endomorphisms $\lambda_u$ corresponding to unitaries $u$ in the normalizer
of the diagonal. Thanks to Power's work \cite{Pow} the structure of this normalizer
is well understood. In the case of restricted Weyl group everything boils down
to analysis of endomorphisms corresponding to permutation unitaries in $M_{n^k}({\mathbb C})$.
Thus one might hope that some straightforward combinatorial manipulations (perhaps
computer aided) with permutations will bring a solution. Unfortunately, as in level
$k$ there are $n^k!$ such permutations, the size of the problem grows too rapidly and
already for very small parameters exceeds computational capacity of modern computers.
For example, $4^3!>10^{89}$ is greater than the number of atoms in the universe.

In order to address both problems mentioned above, we develop a novel combinatorial
approach to the study of permutation related endomorphisms (see Corollary \ref{charaut},
below).  Its essence is reduction of determining invertibility of $\lambda_u$ to
a sequential process involving several steps. In this process labeled rooted trees
are associated to permutations, and certain partial orders on pairs of labels are
considered. These labeled trees also serve as invariants of outer automorphism
classes. Our approach reduces the computational complexity so dramatically as to
allow for relatively simple solution  in such cases as e.g. $n=2$, $k=4$.
Despite $2^4!>2\cdot 10^{13}$ permutations to be considered in this case,
a pen and paper calculation (later verified on a computer) was possible and led
to a complete classification of all permutation related automorphisms of $\O_2$
corresponding to level $4$ (see Subsection \ref{ptwofour}, below).

In examples illustrating our theory we pay particular attention to the case of $\O_2$,
since in some sense this case is the most untractable. Indeed,
concrete examples of permutation related outer automorphisms of $\O_n$, $n\geq3$,
have been known already. It was shown recently in \cite{Sz} that such automorphisms
corresponding to level $2$ generate in ${\rm Out}(\O_n)$ a group containing
free product ${\mathbb Z}_3*{\mathbb Z}_2$. On the other hand, precious little
has been known until now about permutation related automorphisms of $\O_2$. To
the best of our knowledge, the only known example of an outer automorphism of
$\O_2$ of this type was the Archbold's flip-flop \cite{Arc}. Our results explain
why this was so. Namely, new outer automorphisms of $\O_2$ appear only in level $4$,
and to find them one has to sieve through more than $2\cdot10^{13}$ permutations.

Our paper is organized as follows. In Section \ref{start}, we set up notation
and present basic structural results about diagonal preserving automorphisms
of $\O_n$ which follow more or less directly from the works of Cuntz and Power.
In Section \ref{hunt}, we give a general criterion of invertibility of localized
endomorphisms. We also present a criterion for a localized endomorphism to restrict
to an automorphism of the diagonal. In Section \ref{trees}, we develop a labeled tree
approach to the search for permutation related automorphisms. We also discuss
the effect of inner automorphisms and thus show that unlabeled trees are
inner equivalence invariants. In Section \ref{otwo}, we apply the above mentioned
techniques to the case of $\O_2$. In Section \ref{computables}, we present a
more direct approach to finding automorphisms, based on solving certain polynomial matrix
equations. Even though these equations are relatively easy to derive, finding a complete
set of solutions is a highly non-trivial task. We also give tables summarizing
the results of our automorphism search for small values of parameters $n$ and $k$.
These tables were produced through massive computer calculations involving all
of the techniques developed in the present paper.

\vspace{3mm}\noindent{\bf Acknowledgements.}
We owe a great debt of gratitude to Dr Jason Kimberley of Newcastle for his invaluable support
with all the computer calculations, which were performed on MAGMA software.

\section{Setup}\label{start}

If $n$ is an integer greater than 1, then the Cuntz algebra $\O_n$ is a unital, simple
$C^*$-algebra generated by $n$ isometries $S_1, \ldots, S_n$, satisfying
$\sum_{i=1}^n S_i S_i^* = I$ \cite{Cun1}.
We denote by $W_n^k$ the set of $k$-tuples $\alpha = (\alpha_1,\ldots,\alpha_k)$
with $\alpha_m \in \{1,\ldots,n\}$, and by $W_n$ the union $\cup_{k=0}^\infty W_n^k$,
where $W_n^0 = \{0\}$.
We call elements of $W_n$ multi-indices.
If $\alpha = (\alpha_1,\ldots,\alpha_k) \in W_n$, then
$S_\alpha = S_{\alpha_1} \ldots S_{\alpha_k}$ ($S_0 = I$ by convention).
Every word in $\{S_i, S_i^* \ | \ i = 1,\ldots,n\}$ can be uniquely expressed as
$S_\alpha S_\beta^*$, for $\alpha, \beta \in W_n$ \cite[Lemma 1.3]{Cun1}.
If $\alpha \in W_n^k$, then $l(\alpha) = k$, the length of $\alpha$.

$\F_n^k$ is the $C^*$-algebra generated by all words of the form
$S_\alpha S_\beta^*$, $\alpha, \beta \in W_n^k$, and it is isomorphic to the
matrix algebra $M_{n^k}({\mathbb C})$. $\F_n$, the norm closure of
$\cup_{k=0}^\infty \F_n^k$, is the UHF-algebra of type $n^\infty$,
called the core $UHF$ -subalgebra of $\O_n$ \cite{Cun1}.
There exists a faithful conditional expectation $F_0: \O_n \to \F_n$ \cite{Cun1}.

$\D_n$ denotes the diagonal subalgebra of $\O_n$, i.e. the $C^*$-subalgebra
generated by the projections $P_\alpha=S_\alpha S_\alpha^*$, $\alpha \in W_n$.
As remarked by Cuntz \cite{Cun2,CuKr}, $\D_n$ is a maximal abelian subalgebra,
regular both in $\F_n$ and $\O_n$.
$\D_n$ is naturally isomorphic to ${\mathbb C}(X_n)$, where the spectrum $X_n$ is the
collection of infinite words in the alphabet $\{1,\ldots,n\}$ \cite{CuKr}. $X_n$ with the
product topology is a Cantor set, i.e. a compact, metrizable, totally disconnected
space with no isolated points.
There exists a faithful conditional expectation from $\F_n$ onto $\D_n$ and whence
from $\O_n$ onto $\D_n$ as well.
We denote $\D_n^k = \D_n \cap \F_n^k$.

Let $P_n^k$ denote the group of permutations of $W_n^k$.
For $\sigma \in P_n^k$ there is a corresponding unitary $u \in \F_n^k$
(we write $u \sim \sigma$),
$u = \sum_{\alpha \in W_n^k} S_{\sigma(\alpha)} S_\alpha^*$.
We denote $\P_n^k = \{u \ | \ \exists \sigma \in P_n^k, u \sim \sigma\}$
and $\P_n = \cup_{k=0}^\infty \P_n^k$.
We have
${\mathcal N}_{\F_n}(\D_n) = \U(\D_n) \cdot \P_n$
where ${\mathcal N}_{\F_n}(\D_n)$ denotes the (unitary) normalizer of $\D_n$ in $\F_n$ and
$\U(\D_n)$ is the unitary group of $\D_n$.

For $B \subseteq A$ algebras, we denote
${\rm Aut}(A,B)=\{\sigma \in {\rm Aut}(A) \ | \ \sigma(B) = B\}$,
${\rm Aut}_B(A)=\{\sigma \in {\rm Aut}(A) \ | \ \sigma|_B = {\rm id}_B\}$,
${\rm Inn}(A)$ the inner automorphisms,
${\rm Out}(A) = {\rm Aut}(A)/{\rm Inn}(A)$,
and $\pi: {\rm Aut}(A) \to {\rm Out}(A)$ the canonical quotient map.

We recall some notations and results from \cite{Cun2}.
${\rm End}(\O_n)$ is a semigroup (with composition) of unital endomorphisms of
$\O_n$. We have a canonical $\varphi \in {\rm End}(\O_n)$,
$\varphi(a) = \sum_{i=1}^n S_i a S_i^*$. There is a map
$\lambda: \U(\O_n) \to {\rm End}(\O_n)$, determined by $\lambda_u(S_i) = u^* S_i$.
$\lambda$ is a semigroup isomorphism if $\U(\O_n)$ is equipped with the
convolution multiplication $u * w = u \lambda_u(w)$.
The inverse of $\lambda$ is the map $\rho \mapsto \sum_{i=1}^n S_i \rho(S_i^*)$.
Furthermore, ${\rm Aut}(\O_n) = \{\lambda_u \ | \ u^* \in \lambda_u(\O_n)\}$\footnote{
In general, it may happen that $\lambda_u$ is an automorphism but $\lambda_{u^*}$ is not.}
and ${\rm Inn}(\O_n) = \{\lambda_u \ | \ u = \varphi(w)w^*, w \in \U(\O_n)\}$.
The map $\U(\O_n)/{\mathbb T}1 \to {\rm Inn}(\O_n)$,
given by $u \mapsto \lambda_{\varphi(u)u^*} = {\rm Ad}(u)$, is a group isomorphism.
We say that $\lambda_u$ is invertible if $\lambda_u \in {\rm Aut}(\O_n)$.
For $E \subseteq \U(\O_n)$ we denote
$\lambda(E)^{-1} = \{\lambda_u \ | \ u \in E\} \cap {\rm Aut}(\O_n)$.

As shown in \cite{Cun2} we have
${\rm Aut}(\O_n,\D_n) = \lambda({\mathcal N}_{\O_n}(\D_n))^{-1}$
and ${\rm Aut}_{\D_n}(\O_n) = \lambda(\U(\D_n))^{-1} \simeq \U(\D_n)$.
More recently, Power determined in \cite{Pow} (see also \cite{HPP,Pow2})
the structure of ${\mathcal N}_{\O_n}(\D_n)$. Namely, every $w \in {\mathcal N}_{\O_n}(\D_n)$
has a unique decomposition as $w = tu$ with $t \in \U(\D_n)$ and $u$ a finite sum of words.
That is, $u$ is a unitary such that $u = \sum_{j=1}^m S_{\alpha_j}S_{\beta_j}^*$
for some $\alpha_j, \beta_j \in W_n$. Clearly, such unitaries form a group, which
we denote $\S_n$, and this group acts on on $\U(\D_n)$ by conjugation.
Thus, Power's result says that ${\mathcal N}_{\O_n}(\D_n)$ has the structure of semi-direct
product $\U(\D_n) \rtimes \S_n$. Therefore, one obtains the following result \cite{Sz,Ma}.
\begin{theorem} \label{thm1}
${\rm Aut}(\O_n,\D_n) \simeq \U(\D_n) \rtimes \lambda(\S_n)^{-1}.$
In particular, $\lambda(\S_n)^{-1}$ is a subgroup of ${\rm Aut}(\O_n,\D_n)$.
\end{theorem}

Turning back to automorphisms which preserve both the diagonal and the UHF subalgebra,
one easily deduces from the above that
${\mathcal N}_{\F_n}(\D_n) = \U(\D_n) \rtimes \P_n$ as $\P_n = \S_n \cap \F_n$.
Consequently, one has the following result \cite{Sz}.
\begin{theorem}\label{thm1b}
${\rm Aut}(\O_n,\D_n) \cap {\rm Aut}(\O_n,\F_n) = \lambda({\mathcal N}_{\F_n}(\D_n))^{-1}
\simeq \lambda(\U(\D_n)) \rtimes \lambda(\P_n)^{-1} \ . $
In particular, $\lambda(\P_n)^{-1}$ is a subgroup of ${\rm Aut}(\O_n,\D_n) \cap {\rm Aut}(\O_n,\F_n)$.
\end{theorem}

\begin{proof}
At first we show that
${\rm Aut}(\O_n,\D_n) \cap {\rm Aut}(\O_n,\F_n) = \lambda({\mathcal N}_{\F_n}(\D_n))^{-1}.$
If $\lambda_w \in {\rm Aut}(\O_n,\D_n) \cap {\rm Aut}(\O_n,\F_n)$
then it follows from \cite[Prop. 1.5, Prop. 1.2(b)]{Cun2} that
$w \in {\mathcal N}_{\O_n}(\D_n) \cap \F_n = {\mathcal N}_{\F_n}(\D_n)$.
On the other hand,
if $w \in {\mathcal N}_{\F_n}(\D_n)
$ and $\lambda_w \in {\rm Aut}(\O_n)$
then $\lambda_w \in {\rm Aut}(\O_n,\D_n)$ and $\lambda_w(\F_n) \subset \F_n$
and the conclusion follows immediately from \cite[Lemma 2]{Sz}.

Let $u \in \P_n$ and let $\lambda_u$ be invertible.
Then $\lambda_u^{-1}$ belongs to ${\rm Aut}(\O_n,\D_n) \cap {\rm Aut}(\O_n,\F_n)$
and thus $\lambda_u^{-1} = \lambda_z$ with $z \in {\mathcal N}_{\F_n}(\D_n)$.
Thus,
by \cite[Lemma 5.4, (i)]{Pow},
there are $v \in \P_n$ and $y \in \U(\D_n)$ such that $z = vy$.
We have ${\rm id} = \lambda_u \lambda_{vy}$ and hence
$1 = u \lambda_u(v)\lambda_u(y)$.
Thus $\P_n \ni u \lambda_u(v) = \lambda_u(y^*) \in \U(\D_n)$.
Therefore $y = 1$ and consequently $\lambda_u^{-1} = \lambda_v$.
It follows that $\lambda(\P_n)^{-1}$ is a subgroup of
${\rm Aut}(\O_n,\D_n) \cap {\rm Aut}(\O_n,\F_n)$.
Clearly, $\lambda(\P_n)^{-1}$ acts on ${\rm Aut}_{\D_n}(\O_n) = \lambda(\U(\D_n))$
by conjugation.

Now, by Theorem \ref{thm1},
$\lambda_w \in {\rm Aut}(\O_n,\D_n) \cap {\rm Aut}(\O_n,\F_n)$
can be uniquely written as a product of two elements
from $\lambda(\S_n)^{-1}$ and $\lambda(\U(\D_n))$,
$\lambda_w = \lambda_u \lambda_{s^*}$, $u \in \S_n$, $s \in \D_n$.
But then $\lambda_w \lambda_{s}(\F_n)
= \F_n$ and
$u \in \S_n \cap \F_n = \P_n$.
\end{proof}

A slightly weaker version of the following lemma was given in \cite{Sz}.
\begin{lemma}\label{l:2.3}
Let $w \in \P_n$. If $\lambda_w \in {\rm Inn}(\O_n)$ then there exists $u \in \P_n$ such that
$w = \varphi(u)u^*$.
Moreover, for $k \geq 2$, if $w \in \P_n^k$ then $u \in \P_n^{k-1}$.
\end{lemma}

\begin{proof}
The proof of the first statement can be found in \cite{Sz}.
Suppose that $w = \varphi(u) u^* \in \F_n^k$ with $u \in \F_n^h$ for some $h$.
Observe that if $h \geq k$ then
$\F_n^h \supset \F_n^k$ so that $\varphi(u) \in \F_n^h$ and $u \in \F_n^{h-1}$.
Therefore $h < k$ and necessarily one must have $h = k-1$.
\end{proof}

Since
$\P_n \simeq \lambda(\P_n)^{-1} \cap {\rm Inn}(\O_n)$ via $u \mapsto {\rm Ad}(u)$ \cite{Cun2},
there exists an exact sequence
\begin{equation}
1 \to \P_n \to \lambda(\P_n)^{-1} \to \pi(\lambda(\P_n)^{-1}) \to 1 \ .
\end{equation}

The natural inclusion $\P_n^k \subset \P_n^{k+m}$ corresponds to the embedding
$P_n^k \hookrightarrow P_n^{k+m}$, $\phi \mapsto \phi \times {\rm id}_m$,
where ${\rm id}_m$ denotes the identity on $W_n^m$
(we have $W_n^{k+m} = W_n^k \times W_n^m$).
With this identification $P_n = \cup_k P_n^k$ becomes a group isomorphic to $\P_n$.
We note that the imbedding $\P_n^k \hookrightarrow \P_n^{k+m}$,
$u \mapsto \varphi^m(u)$, corresponds to the imbedding
$P_n^k \hookrightarrow P_n^{k+m}$, $\phi \mapsto {\rm id}_m \times \phi$.
If $\phi \in P_n^k$ and $r \geq 1$ then we define $\phi^{(r)} \in P_n^{k+r-1}$ as
\begin{equation}
\phi^{(r)}
= (id_{r-1} \times \phi) \, ({\rm id}_{r-2} \times \phi \times {\rm id}_1)
\ldots (\phi \times {\rm id}_{r-1}) \ .
\end{equation}
In particular, $\phi^{(1)} = \phi$.
For $u \in \P_n^k$, $u \sim \phi$, $w = \varphi(u)u^*$
(i.e., $\lambda_w = {\rm Ad}(u)$), $w \sim \psi$,
we have $\psi = ({\rm id}_1 \times \phi)(\phi^{-1} \times {\rm id}_1)$,
$\psi^{(r)} = ({\rm id}_r \times \phi) (\phi^{-1} \times {\rm id}_r)$
and thus $\psi^{(k)} = \phi^{-1} \times \phi$.

The map $\P_n^k \times \P_n^r \to \P_n^{k+r-1}$,
$(u,w) \mapsto u * w = u\lambda_u(w)$ corresponds to the map
$P_n^k \times P_n^r \to P_n^{k+r-1}$,
\begin{equation}
(\alpha,\beta) \mapsto \alpha * \beta
= (\alpha \times {\rm id}_{r-1})
(\alpha^{(r)})^{-1} (\beta \times {\rm id}_{k-1}) \alpha^{(r)} \ .
\end{equation}
If a permutation $\phi \in P$ is $*$-invertible, then we denote its inverse
by $\overline\phi$.

\medskip

For later use we highlight a simple but suggestive
reformulation of the second statement in Theorem \ref{thm1b}
with a slight improvement provided by Corollary \ref{bound} below.

\begin{proposition}
\label{localizedinverse}
Let $w \in \P_n^k$ and suppose that $\lambda_w \in {\rm Aut}(\O_n)$,
then
$\lambda_w^{-1}$ is also induced by a unitary in
$\P_n^h$, with $h \leq n^{2(k-1)}$.
\end{proposition}

Following \cite{CP}, 
an endomorphism of the Cuntz algebra $\O_n$
of the form $\lambda_u$ with $u$ a unitary in $\cup_k \F_n^k$ is called ``localized''.
Of course, by the very definition all the $\lambda_w$'s with $w \in \P_n^k$
are localized endomorphisms and,
by the above,
automorphisms induced by permutation unitaries are examples of
localized automorphisms whose inverse is
(induced by a permutation unitary and thus)
still localized.

\medskip

Now the natural question arises whether one can
find an effective algorithm to identify
all the permutation unitaries inducing automorphisms of the Cuntz algebra.
Eventually, one might also like to give a closer look at the properties of these
(possibly outer) automorphisms.
Especially, one interesting problem is to determine the structure of the groups
$\pi(\lambda(\P_n)^{-1})$ for $n \geq 2$.
As
shown in \cite[Example 9]{Sz},
 the groups $\pi(\lambda(\P_n)^{-1})$ for $n \geq 3$ are quite ``big''
in the sense that they contain non-amenable subgroups,
notably ${\mathbb Z}_2 * {\mathbb Z}_3$.
The same question for $n=2$ is more subtle.
On the basis of general results \cite{Ro,i2,Matui}, it has been known for some time
that the automorphism group of $\O_2$ is in some sense considerably ``smaller''.
Our computations provide a very concrete evidence to this effect.

\section{Searching for automorphisms}\label{hunt}

\subsection{Invertibility of localized endomorphisms}

Let $w \in \P_n^k$ be a permutation unitary. We set
\begin{equation}
B_w = \{w,\varphi(w),\ldots,\varphi^{k-2}(w)\}' \cap \F_n^{k-1}
\end{equation}
if $k \geq 2$ and $B_w = {\mathbb C}1$ otherwise.
That is, $b \in \F_n^{k-1}$ is in $B_w$ if and only if,
for any $\alpha, \beta \in W_n^l$, $l \in \{0,\ldots,k-2\}$,
$S_\alpha^* b S_\beta$ commutes with $w$.
Of course, $B_w (= B_{w^*})$ is a unital $*$-subalgebra of $\F_n^{k-1}$.
Notice that if $b \in B_w$ then $\lambda_w(b) = b$.
Also, for $i,j \in \{1,\ldots,n\}$ we define maps
$a_{ij}^w: \F_n^{k-1} \to \F_n^{k-1}$ by
\begin{equation}
a_{ij}^w(x) = S_i^* w x w^* S_j, \;\;\; x \in \F_n^{k-1}.
\end{equation}
We denote $V_w = \F_n^{k-1}/B_w$.
Since $a_{ij}^w(B_w) \subseteq B_w$, there are induced maps
$\tilde{a}_{ij}^w: V_w \to V_w$.
We define $A_w$ as the subring of $\L(V_w)$ generated by
$\{\tilde{a}_{ij}^w \ | \ i,j = 1,\ldots, n\}$.

In the sequel we elaborate on \cite[Theorem 7]{Sz}
and provide further insight on that matter discussing a powerful extension of the argument.
\begin{lemma}\label{nil}
If $w \in \P_n$ then $\lambda_w$ is invertible if and only if $A_w$ is nilpotent.
\end{lemma}

\begin{proof}
Necessity. Let $w\in \P_n^k$ and $\lambda_w$ be invertible.
By Proposition \ref{localizedinverse}, $\lambda_w^{-1}$ is then induced by some
(permutation) unitary in some finite matrix algebra.
Let $\lambda_w^{-1}(\F_n^{k-1}) \subseteq \F_n^l$.
For $a \in \F_n^l$ the sequence
${\rm Ad}(w^*\varphi(w^*)\ldots\varphi^m(w^*))(a)$ stabilizes from $m = l-1$
at $\lambda_w(a)$. Consequently, for any $b \in \F_n^{k-1}$ the sequence
${\rm Ad}(\varphi^m(v) \ldots \varphi(w)w)(b)$ stabilizes from $m = l-1$ at
$\lambda_w^{-1}(b)$.
There are $c_{\gamma\rho} \in {\mathbb C}1$ such that
$$\sum_{\gamma,\rho \in W_n^l} S_\gamma c_{\gamma\rho}(b) S_\rho^*
= {\rm Ad}(\varphi^{l-1}(w) \ldots \varphi(w)w)(b) \in \F_n^l \ . $$
If $\alpha = (i_1,\ldots,i_l)$, $\beta = (j_1,\ldots,j_l)$,
$T_{\alpha,\beta} = a_{i_l j_l}^w \ldots a_{i_1 j_1}^w$,
and $b \in \F_n^{k-1}$, then $T_{\alpha,\beta}(b) = c_{\alpha\beta}(b)
\in {\mathbb C}1 \subset B_w$. Consequently, $A_w^l = 0$.

\medskip
Sufficiency.
Let $w \in \P_n^k$ and assume that $A_w^l = 0$.
Let $b \in \F_n^{k-1}$ and $T_{\alpha,\beta}$ as above.
By hypothesis, $T_{\alpha,\beta}(b)$
commutes with $\varphi^m(w)$ for any $m$.
Hence, if $r \geq 1$, we have
\begin{align*}
{\rm Ad}(\varphi^{l-1+r}(w) \ldots \varphi(w)w) (b)
& =
{\rm Ad}(\varphi^{l-1+r}(w) \ldots \varphi^l(w) )
\Big(\sum_{\alpha,\beta \in W_n^l} S_\alpha T_{\alpha,\beta}(b) S_\beta^* \Big) \\
& =
\sum_{\alpha,\beta \in W_n^l} S_\alpha
{\rm Ad}(\varphi^{r-1}(w) \ldots w) (T_{\alpha,\beta}(b)) S_\beta^* \\
& = \sum_{\alpha,\beta \in W_n^l} S_\alpha T_{\alpha,\beta}(b) S_\beta^* \ .
\end{align*}
Thus, for any $b \in \F_n^{k-1}$, the sequence
${\rm Ad}(\varphi^m(w) \ldots \varphi(w)w)(b)$ stabilizes from $m = l-1$.
Let $w^* = \sum_{i,j = 1}^n S_i b_{ij} S_j^*$, $b_{ij} \in \F_n^{k-1}$.
By the above, the sequence
\begin{align*}
{\rm Ad}(\varphi^m(w) \ldots \varphi(w)w) (w^*)
& = \sum_{i,j = 1}^n {\rm Ad}(\varphi(\varphi^{m-1} \ldots \varphi(w)w))
(S_i b_{ij}S_j^*) \\
& = \sum_{i,j} S_i {\rm Ad}(\varphi^{m-1}(w) \ldots \varphi(w)w)(b_{ij}) S_j^*
\end{align*}
stabilizes from $m=l$ at $\lambda_w^{-1}(w^*)$ and hence $\lambda_w$ is invertible.
\end{proof}

In turn,
inspection of the proof shows that a similar characterization holds true
for any unitary $u \in \F_n^k$ such that $\lambda_u$ is invertible
with localized inverse.
(If $\lambda_u^{-1} = \lambda_v$ with $v \in \F_n^h$ one can choose $l = k + h -2$ in the above argument).
Moreover, we can adapt some arguments from section 6 of \cite{CP} to our situation.
We denote by $H$ the linear span of the $S_i$s.
Given a unitary $u \in \F_n^k$,
let us define inductively
\begin{equation}
\Xi_0 = \F_n^{k-1}, \quad
\Xi_{r} = \lambda_u(H)^* \Sigma_{r-1} \lambda_u(H) \, , \ r \geq 1 \ ,
\end{equation}
that is $\Xi_r = (\lambda_u(H)^r)^* \F_n^{k-1} (\lambda_u(H))^r$.
It readily follows that $(\Xi_r)_r$ is nonincreasing sequence of subspaces
of $\F_n^{k-1}$
that stabilizes at the first value $p$ for which $\Xi_p = \Xi_{p+1}$.
Let $\Xi_u := \bigcap_r \Xi_r = \Xi_p$.

\begin{theorem}\label{main}
Let $u$ be a unitary in $\F_n^k$ for some $k \geq 1$.
Then the following conditions are equivalent:
\begin{itemize}

\item[(1)] $\lambda_u$ is invertible with localized inverse;

\item[(2)]
the sequence of unitaries
$$\big({\rm Ad}(\varphi^m(u)\varphi^{m-1}(u) \ldots \varphi(u)u)(u^*)\big)_{m \geq 1}$$
eventually stabilizes;

\item[(3)] $A_u$ is nilpotent;

\item[(4)] $\Xi_u \subseteq B_u$;

\item[(5)] $\Xi_u = {\mathbb C}1$.
\end{itemize}
\end{theorem}

\begin{proof}
(1) $\Rightarrow$ (2):
let $v \in \F_n^h$ be such that $\lambda_u \lambda_v = {\rm id}$.
Thus $u \lambda_u(v) = 1$, that is
$$u^* \varphi(u^*) \cdots \varphi^{m}(u^*) v \varphi^{m}(u) \cdots \varphi(u)u = u^*$$
for every $m \geq h-1$.\\
(2) $\Rightarrow$ (1):
Suppose that there exists some positive integer $l$ for which it holds
$$\varphi^m(u) \ldots \varphi(u) u^* \varphi(u^*) \ldots \varphi^m(u^*)
= \varphi^l(u) \ldots \varphi(u) u^* \varphi(u^*) \ldots \varphi^l(u^*)$$
for every $m \geq l$.
Call $v$ the resulting unitary, clearly in $\F_n^{k+l}$.
Then $u \lambda_u(v) = u (u^* \ldots \varphi^{k+l-1}(u^*)) v (\varphi^{k+l-1}(u) \ldots u)
= u u^* = 1$, and therefore $\lambda_v = \lambda_u^{-1}$.

The equivalence of (1) and (3) follows by Lemma \ref{nil}, mutatis mutandis.
As (4) is nothing but a reformulation of the nilpotency condition, (3) and (4) are clearly equivalent.
\\
(4) $\Rightarrow$ (5): suppose that $\Xi_u \subset B_u$. Then $u \Xi_u u^* = \Xi_u$.
Now, as recalled above, one has
\begin{align*}
\Xi_u & = (H^p)^* \varphi^{p-1}(u) \ldots u \F_n^{k-1} u^* \ldots \varphi^{p-1}(u^*) H^p \\
& = (H^{p+1})^* \varphi^{p}(u) \ldots u \F_n^{k-1} u^* \ldots \varphi^{p}(u^*) H^{p+1} ,
\end{align*}
and by assumption $\Xi_u$ is also equal to
$(H^p)^* \varphi^{p}(u) \ldots u \F_n^{k-1} u \ldots \varphi^{p}(u) H^p$.
It readily follows that $\Xi_u = H^* \Xi_u H$ and thus $\Xi_u = {\mathbb C}.$ \\
(5) $\Rightarrow$ (4): obvious.
\end{proof}
Note that implication (1) $\Rightarrow (5)$ in the above theorem also follows by \cite[Proposition 6.1]{CP},
where we take as $\Phi$ the (normal extension of the) localized automorphism ${\lambda_u}^{-1}$.\footnote{
We warn the reader about a slightly confusing change in the conventions.
The $\lambda_u$ in \cite{CP} corresponds to $\lambda_{u^*}$ here.}

\begin{corollary}\label{bound}
Let $u \in \F_n^k$ be a unitary satisfying the equivalent conditions of
Theorem \ref{main}.
Then $\lambda_u^{-1}$ is induced by a unitary $v \in \F_n^h$ with $h = n^{2(k-1)}$.
\end{corollary}

\begin{proof}
As the sequence of finite dimensional subspaces
$\F_n^{k-1} \supset K^* \F_n^{k-1} K \supset {K^*}^2 \F_n^{k-1} K^2 \supset \ldots $
is decreasing until it stabilizes to $\mathbb C$,
$\dim(\F_n^{k-1}) = n^{2(k-1)}$ and at each step the dimension drops by one at least,
one has $(K^*)^p \F_n^{k-1} K^p = {\mathbb C}$ for some $p \leq n^{2(k-1)}-1$.

Next observe that
$${K^*}^{p+1} \F_n^k K^{p+1}
= {K^*}^p \F_n^{k-1} K^p = {\mathbb C} \ . $$
That is,
$${H^*}^{p+1} \varphi^p(u) \cdots u \F_n^{k} u^* \cdots \varphi^p(u^*) H^{p+1} = {\mathbb C}$$
and
$v:=\varphi^p(u) \cdots u u^* u^* \cdots \varphi^p(u^*) \in \F_n^{p+1}$.
This shows the statement.
\end{proof}

At this stage it is not clear whether it is possible to improve the exponential
bound on $h$ in the last corollary.
This would be rather useful for computational purposes.

\subsection{Automorphisms of the diagonal}

It follows from \cite[Proposition 1.5]{Cun2} that if $w \in \F_n^k$
is in the normalizer of the diagonal subalgebra $\D_n$ then for
$\lambda_w$ to be invertible it is necessary that $\lambda_w(\D_n) = \D_n$.
It turns out that the method of Lemma \ref{nil} and Theorem \ref{main}
can also provide a criterion of invertibility of the restriction
of such an endomorphism $\lambda_w$ to the diagonal $\D_n$.

Indeed, let $w\in\F_n^k\cap{\mathcal N}_{\O_n}(\D_n)$. Then both $\D_n^{k-1}$
and $B_w\cap\D_n^{k-1}$ are invariant subspaces for all the operators
$a_{ij}^w$ associated with $w$. Denote the restriction of $a_{ij}^w$ to
$\D_n^{k-1}$ by $b_{ij}^w$. Each $b_{ij}^w$ induces a linear
transformation $\tilde{b}_{ij}^w:V_w^D\rightarrow V^D_w$, where
$V^D_w=\D_n^{k-1}/B_w\cap\D_n^{k-1}$. We denote by $A^D_w$ the
subring of ${\mathcal L}(V^D_w)$ generated by $\{\tilde{b}_{ij}^w \ | \ i,j=1,\ldots,n\}$.
Also, we consider the subspace of $\D_n^{k-1}$ defined by
$\Xi_w^D := \bigcap_r (K^*)^r \D_n^{k-1}K^r$,
where $K$ is
the linear span of $w^* S_1, \ldots , w^* S_n$.

\begin{theorem}\label{nildiag}
Let $w$ be a unitary in $\F_n^k \cap {\mathcal N}_{\O_n}(\D_n)$.
If the ring $A_w^D$ is nilpotent then $\lambda_w$ restricts to an automorphism of $\D_n$.
More precisely, the following conditions are equivalent:
\begin{itemize}
\item[(1)] $\lambda_w$ restricts to an automorphism of the algebraic
part $\cup_s \D_n^s$ of $\D_n$;
\item[(2)] the ring $A_w^D$ is nilpotent;
\item[(3)] $\Xi_w^D \subseteq B_w \cap \D$;
\item[(4)] $\Xi_w^D = {\mathbb C}1$.
\end{itemize}
\end{theorem}
\begin{proof}
We only give details of the proof of implication (2) $\Rightarrow$ (1).
The other implications are established through arguments very similar to those
of Lemma \ref{nil} and Theorem \ref{main}.

Suppose that $A_w^D$ is nilpotent.
We show by induction on $r\geq k$ that all $\D_n^r$ are in
the range of $\lambda_w$
restricted to $\cup_s \D_n^s$.

If $x\in\D_n^k$ then the same argument as in the proof of sufficiency part
in Lemma \ref{nil} shows that $x$ belongs to
$\lambda_w(\cup_s \D_n^s)$. In fact,
the sequence ${\rm Ad}(\varphi^m(w) \ldots \varphi(w)w)(x)$ stabilizes at
$\lambda_w^{-1}(x) \in \cup_s \D_n^s$.

For the inductive step, suppose that $r\geq k$ and $\D_n^r\subset\lambda_w( \cup_s \D_n^s)$.
Since $\D_n^{r+1}$ is generated by $\D_n^r$ and $\varphi^r(\D_n^1)$, it suffices
to show that $\varphi^r(y)$ belongs to $\lambda_w(\cup_s \D_n^s)$ for all $y\in\D_n^1$.
However, $\varphi^r(y)$ commutes with $w$ and $\varphi^{r-1}(y)\in\D_n^r$ is in
$\lambda_w (\cup_s \D_n^s)$. Thus, we see that the sequence
$$ {\rm Ad}(\varphi^m(w)\ldots \varphi(w)w)(\varphi^r(y))=\varphi({\rm Ad}(\varphi^{m-1}(w)
   \ldots\varphi(w)w)(\varphi^{r-1}(y))) $$
stabilizes at $\lambda_w^{-1}(\varphi^r(y)) \in \cup_s \D_n^s$.
\end{proof}

\section{Applications of labeled trees to the search for automorphisms}\label{trees}

Let $w \in \P_n^k$.
Take $\{S_\alpha S_\beta^*\}_{\alpha,\beta \in W_n^{k-1}}$,
a basis of $\F_n^{k-1}$, so that $\{S_\alpha S_\alpha^*\}$
are the first block of the basis.
With respect to this basis, each $a_{ij}^w$, $i,j \in \{1,\ldots,n\}$ has a matrix
\begin{equation}
a_{ij}^w = \begin{pmatrix} b_{ij}^w & c_{ij}^w \\ 0 & d_{ij}^w \end{pmatrix}
\end{equation}
with entries in $\{0,1\}$,
as $a_{ij}^w(S_\alpha S_\beta^*) = \sum_m S_i^* S_{\sigma(\alpha,m)} S_{\sigma(\beta,m)}^* S_j$,
where $w \sim \sigma$.

\medskip

In the sequel of this section, we will explain how the condition that $\lambda_w
\in {\rm Aut}(\O_n)$ for $w \in \P_n^k$ translates in terms of the $a_{ij}^w$'s.
In turn, this boils down to two separate arguments for the (sub-)matrices $[b]$
and $[d]$.
As a matter of fact, $[c]$ turns out to be irrelevant for the following
discussion.
Indeed,
since $a_{ij}^w(I) \in \{I,0\}$, each $a_{ij}^w$ gives rise to a map from $\F_n^{k-1}/{\mathbb C}$
to itself, whose matrix has a block form
\begin{equation}
\begin{pmatrix} \hat{b}_{ij}^w & * \\ 0 & d_{ij}^w \end{pmatrix} \ .
\end{equation}
It is an immediate corollary of Lemma \ref{nil} and Theorem \ref{main} that $\lambda_w$
is invertible if and only if both rings generated by $\{\hat{b}_{ij}^w \ | i,j = 1,\ldots, n\}$
and by $\{d_{ij}^w \ | \ i,j = 1,\ldots, n \}$, respectively, are nilpotent.
Furthermore, it follows from Proposition \ref{nildiag} that nilpotency
of the ring generated by $\{\hat{b}_{ij}^w  \ | \, i,j = 1,\ldots, n\}$ implies that
endomorphism $\lambda_w$ restricts to an automorphism of $\D_n$.

\subsection{Upper left corner $[b]$}

The plan of this subsection is as follows.
We first convert the matrix $[b]$ into functions on indices.
Trees then pop up as diagrams of these functions.
Next we discuss labeling.
The automorphism condition will lead us to trees with a suitable labeling,
that is inducing a certain partial order relation.

\medskip

If $i \neq j$ then $b_{ij}^w = 0$. Hence we can write $b_i^w := b_{ii}^w$. Since
\begin{equation}
b_i^w (S_\alpha S_\alpha^*) = \sum_m S_i^* w S_\alpha S_m S_m^* S_\alpha^* w^*
S_i \ ,
\end{equation}
$w S_\alpha S_m S_m^* S_\alpha^* w^*$ being a minimal projection in $\D_n^k$,
it follows that each column of $b_i^w$ has at most $n$ non-zero entries
but fixing a column and summing over $i$ we get exactly $n$. Furthermore,
since $b_i^w(1) = 1$, we have
\begin{equation}
\sum_\alpha S_\alpha S_\alpha^* = \sum_\alpha b_i^w(S_\alpha S_\alpha^*)
\end{equation}
and hence each row of $b_i^w$ has exactly one 1 and the rest 0.

\medskip

Suppose that $\lambda_w$ is an automorphism of $\O_n$. Then
equivalence of conditions (1), (3) and (5) of Theorem \ref{main}
easily implies the following condition
on the left-upper corner of the matrix $[a_{ij}^w]$:
sufficiently long products of the operators $\{b_i^w \ | \ i=1,\ldots,n\}$
have the form
\begin{equation}\label{constb}
\begin{pmatrix}
\lambda_1 & \lambda_2 & \cdots & \lambda_{n^{k-1}} \\
\lambda_1 & \lambda_2 & \cdots & \lambda_{n^{k-1}} \\
\vdots & \vdots &  & \vdots \\
\lambda_1 & \lambda_2 & \cdots & \lambda_{n^{k-1}} \\
\end{pmatrix} \ ,
\end{equation}
that is they are constant along the columns. However, since for any $i$
each row of $b_i^w$ contains exactly one non-zero entry, the same is true for
products of $\{b_i^w\}$s.
Thus each of the above matrices as in (\ref{constb}) must actually have the form
\begin{equation}\label{const01}
\begin{pmatrix}
0 & \cdots & 1 & & \cdots & 0 \\
0 & \cdots & 1 & & \cdots & 0 \\
\vdots & & \vdots & & & \vdots \\
0 & \cdots & 1 & & \cdots & 0 \\
\end{pmatrix} \ ,
\end{equation}
i.e. one column of 1's and 0's elsewhere.

Since each row of the matrix $b_i^w$  has 1 exactly in one
column and 0's elsewhere, the $b_i^w$ can be identified with a function $f_i^w:
W_n^{k-1} \to W_n^{k-1}$ defined by
\begin{equation}\label{ef}
f_i^w(\alpha) = \beta
\end{equation}
whenever $b_i^w$ has 1 in $\alpha$-$\beta$ entry.
Suppose  that  $w$ comes from a permutation $\sigma$. Then
\begin{eqnarray}
f_i^w(\alpha) = \beta & \Longleftrightarrow & \exists
m \mbox{ such that } (i,\alpha)=\sigma(\beta,m) \\
    & \Longleftrightarrow & S_\alpha S_\alpha^* \leq S_i^* w S_\beta S_\beta^* w^* S_i \ .
\end{eqnarray}
It is not difficult to verify that the product $b_i^w b_j^w$ corresponds to the composition
$f_j^w \circ f_i^w$ (in reversed order of $i$ and $j$).
In what follows we often omit superscript $w$ in $f_i^w$ when no confusion may arise.

We omit an easy proof of the following lemma.
\begin{lemma}
The ring generated by $\{\hat{b}_i^w \ | \ i=1,\ldots,n\}$
is nilpotent if and only if all sufficiently long composition products of
mappings $\{f_i\;|\;i=1,\ldots,n\}$ have ranges consisting of a single element.
\end{lemma}
\begin{lemma}\label{necfi}
A necessary condition of nilpotency of the ring
generated by $\{\hat{b}_i^w \ | \ i=1,\ldots,n\}$
is that each $f_i$ must have the following structure:
\begin{itemize}
\item exactly one fixed-point;
\item no periodic orbits of length $\geq 2$.
\end{itemize}
\end{lemma}
\begin{proof}
The first condition clearly follows by considering, for any given index $i$,
only powers of the matrix $b_i$ or, equivalently,
compositions of the same function $f_i$.
The second condition follows since otherwise some power of $b_i$ would have more
than one fixed-point.
\end{proof}

From this lemma we deduce that the diagrams of the $f_i$'s are rooted trees,
where the root corresponds to the unique fixed point.
By diagram we mean the graph with vertices labeled by elements of $W_n^{k-1}$
and with a directed edge from vertex $\alpha$ to vertex $\beta$ if
$f_i(\alpha) = \beta$. By convention, we do not include in the diagram
the loop from the root (fixed point) to itself.

\begin{example}\label{idtree}
{\rm The pair of labeled trees corresponding to
$\sigma = {\rm id}$ in $P_2^3$. All the edges are downward oriented. }

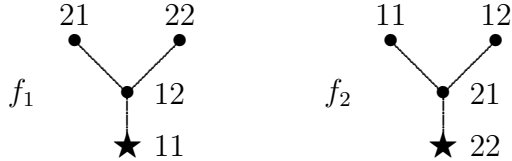
\[ \beginpicture
\setcoordinatesystem units <0.7cm,0.7cm>
\setplotarea x from 4 to 9, y from -1 to 1

\put {$\bullet$} at -1 1
\put {$\bullet$} at 1 1
\put {$\bullet$} at 0 0
\put {$\bigstar$} at 0 -1

\setlinear
\plot -1 1 0 0 /
\plot 1 1 0 0 /
\plot 0 0 0 -1 /
\put {$f_1$} at -2 0
\put {$21$} at -1 1.5
\put {$22$} at 1 1.5
\put {$12$} at 0.8 0
\put {$11$} at 0.8 -1

\put {$\bullet$} at 5 1
\put {$\bullet$} at 7 1
\put {$\bullet$} at 6 0
\put {$\bigstar$} at 6 -1

\setlinear
\plot 5 1 6 0 /
\plot 7 1 6 0 /
\plot 6 0 6 -1 /
\put {$f_2$} at 4 0
\put {$11$} at 5 1.5
\put {$12$} at 7 1.5
\put {$21$} at 6.8 0
\put {$22$} at 6.8 -1

\endpicture \]

\end{example}

\begin{example}\label{bogolubov}
{\rm Let $u\in{\mathcal P}_n^1$, so that $\lambda_u$ is a Bogolubov
automorphism of $\O_n$. If we view $u$ as an element of ${\mathcal P}_n^k$
then all $n$ unlabeled trees corresponding to $u$ are identical; the root
receives $n-1$ edges from other vertices, each other vertex receives either
none or $n$ edges, and the height of the tree (the length of the longest path
ending at the root) is minimal and equal to $k-1$. In particular, all
such unitaries have the corresponding $n$-tuples of unlabeled trees
identical with those of the identity.}
\end{example}

\medskip

\begin{lemma}\label{nsfi}
The ring generated by $\{\hat{b}_i^w \ | \ i=1,\ldots,n\}$ is nilpotent if
and only if there exists a partial order $\leq$ on the cartesian product
$W_n^{k-1} \times W_n^{k-1}$
such that:
\begin{itemize}
\item[(i)] Each element of the diagonal $(\alpha,\alpha)$ is minimal;
\item[(ii)] Each $(\alpha,\beta)$ is bounded below by some diagonal element;
\item[(iii)] For every $i$ and all $(\alpha,\beta)$ such that $\alpha \neq \beta$, we have
\begin{equation}
(f_i(\alpha),f_i(\beta)) \leq (\alpha,\beta) \ .
\end{equation}
\end{itemize}
\end{lemma}
\begin{proof}
Suppose that the ring generated by $\{\hat{b}_i^w \ | \ i=1,\ldots,n\}$ is nilpotent.
Define a relation $\leq$ as follows. For any $\alpha$, $(\alpha,\alpha)\leq(\alpha,
\alpha)$. If $\gamma \neq \delta$ then
$(\alpha,\beta) \leq (\gamma,\delta) $ if and only if there exists a sequence
$j_1,\ldots,j_d$, possibly empty, such that $\alpha = f_{j_1} \circ \cdots \circ
f_{j_d}(\gamma)$ and $\beta = f_{j_1} \circ \cdots \circ f_{j_d}(\delta)$.

Reflexivity and transitivity of $\leq$ are obvious.
Suppose $(\alpha,\beta) \leq (\gamma,\delta)$ and $(\gamma,\delta) \leq (\alpha,\beta)$.
If $(\alpha,\beta)\neq(\gamma,\delta)$ then, by definition of $\leq$,
$\alpha \neq \beta$, $\gamma \neq \delta$ and there exist indices
$j_1,\ldots,j_d$, $k_1,\ldots,k_h$ such that
$(\alpha,\beta) = (f_{j_1} \circ \cdots \circ f_{j_d})(\gamma,\delta)$ and
$(\gamma,\delta) = (g_{k_1} \circ \cdots \circ g_{k_h})(\alpha,\beta)$.
Then $(\alpha,\beta) = (f_{j_1} \circ \cdots \circ f_{j_d} \circ g_{k_1} \circ \cdots
\circ g_{k_h})(\alpha,\beta)$. That is, $t = f_{j_1} \circ \cdots
\circ f_{j_d} \circ g_{k_1} \circ \cdots \circ g_{k_h}$ has
two distinct fixed points,  a contradiction. Thus $(\alpha,\beta) = (\gamma,\delta)$
and $\leq$ is also antisymmetric.

We must still show that each $(\alpha,\beta)$, $\alpha \neq \beta$,
is bounded below by a diagonal
element. If not, then counting shows that there exists a sequence
$f_1, \ldots, f_d$ such that $(\alpha,\beta)
= f_1 \circ \ldots \circ f_d(\alpha,\beta)$ and again, $f_1 \circ \ldots \circ f_d$ has two
distinct fixed points.

\smallskip

Conversely, suppose such a partial order exists. We must show that each
composition of sufficiently many functions $\{f_i\}$ has range consisting of exactly one
element. By counting, to this end it suffices to show that for
any subset $X \subseteq W_n^{k-1}$
with at least two elements and a sufficiently large $r$
the set $f_1 \circ \cdots \circ f_r(X)$ has at least one element less
than $X$. To see this take any two distinct elements $\alpha \neq \beta \in X$. Then, by
the conditions on $\leq$, eventually $f_1 \circ \cdots \circ f_r (\alpha)=
f_1 \circ \cdots \circ f_r(\beta)$, and this does the job.
\end{proof}

Proposition \ref{nildiag} and Lemma \ref{nsfi} yield the following.
\begin{corollary}
Let $w\in\P_n^k$.
If there exists a partial order on $W_n^{k-1} \times W_n^{k-1}$
satisfying conditions of Lemma \ref{nsfi} then endomorphism $\lambda_w$
restricts to an automorphism of $\D_n$.
\end{corollary}

The relation used in Lemma \ref{nsfi} can be explicitly described as follows.
We have that $(\alpha,\beta) \leq (\gamma,\delta)$ if and only if
either $\alpha=\gamma$ and $\beta=\delta$, or $\gamma \neq \delta$  and
there exist $i_0,\ldots,i_r$ such that
\begin{equation}\label{rel}
\begin{array}{cc}
(i_0,\gamma) = \sigma(\gamma_1,k_1), & (i_0,\delta) = \sigma(\delta_1,h_1) \\
(i_1,\gamma_1) = \sigma(\gamma_2,k_2), & (i_1,\delta_1) = \sigma(\delta_2,h_2) \\
\ldots & \ldots  \\
(i_r,\gamma_r) = \sigma(\alpha,k_{r+1}), & (i_r,\delta_r) = \sigma(\beta,h_{r+1}) .
\end{array}
\end{equation}

\medskip

In order to give an equivalent reformulation of Lemma \ref{nsfi} we define inductively
a nested sequence of subsets $\Sigma_m^w$ of $W_n^{k-1}\times W_n^{k-1}$, as follows.
\begin{align}
\Sigma_0^w & = \{(\alpha,\alpha) \ | \ \alpha \in W_n^{k-1}\}, \\
\Sigma_{m+1}^w & = \{(\alpha,\beta) \ | \ (f_i(\alpha),f_i(\beta)) \in \Sigma_m^w ,
    \,i=1,\dots,n\} \cup \Sigma_m^w.
\end{align}
We omit an easy proof of the following proposition.
\begin{proposition}\label{condb}
The relation $\leq$ defined by (\ref{rel}) satisfies conditions of
Lemma \ref{nsfi} if and only if
\begin{equation}\label{bcond}
\bigcup_m \Sigma_m^w = W_n^{k-1} \times W_n^{k-1} \ .
\end{equation}
\end{proposition}

\subsection{Effect of inner automorphisms}\label{inners}

If $w \sim \sigma \in \P_n^k$ ($w = \sum S_{\sigma(\alpha)}S_\alpha^*$) and
$u \sim \phi \in \P_n^{k-1}$ then ${\rm Ad}(u) \lambda_w = \lambda_{\varphi(u)wu^*}$
and $\varphi(u) w u^* \sim (1 \times \phi) \sigma (\phi^{-1} \times 1)$.

Let $f_i$ and $g_i$ be the self-mappings of $W_n^{k-1}$ corresponding to
$w$ and $\varphi(u) w u^*$, respectively, as in (\ref{ef}).
Then $(i,\alpha) = \sigma(\beta,m)$ if and only if
$(i,\phi(\alpha)) = (1 \times \phi) \sigma (\phi^{-1} \times 1)(\phi(\beta),m)$.
Thus $f_i(\alpha)=\beta$ if and only if $g_i(\phi(\alpha)) = \phi(\beta)$.
That is,
\begin{equation}\label{inner}
g_i = \phi f_i \phi^{-1} \ , \quad i = 1, \ldots,n \ .
\end{equation}
Consequently,  the action of inner automorphisms corresponds to permutation of labels. Thus,
combining this observation with Lemma \ref{l:2.3} we obtain the following.
\begin{proposition}\label{outclasses}
Suppose that $u,w\in\P_n^k$ and both $\lambda_u$ and $\lambda_w$ are automorphisms
of $\O_n$. If there exists an $i$ such that the tree corresponding to $f_i^u$ is not
isomorphic to the tree of $f_i^w$ (as directed tree, no labeling involved) then
$\lambda_u$ and $\lambda_w$ give rise to distinct elements of ${\rm Out}(\O_n)$.
\end{proposition}

Now the following question arises: how many distinct permutations $\tau\in P_n^k$ give
rise to the same collection of labeled trees as $\sigma$? The structure and labels
on the trees $\{f_1,\ldots,f_n\}$ corresponding to $\sigma$ are determined by
identity $(i,\alpha) = \sigma(\beta,m)$, in which $m\in\{1,\ldots,n\}$ can
be chosen freely. Thus, simple counting leads to the following.
Given any $w \in \P_n^k$ with corresponding functions $\{f_i^w\}$, there
are exactly $n!^{n^{k-1}}$ elements $u$ of $\P_n^k$ yielding identical maps
$f_i^u=f_i^w$.

\subsection{Lower right corner $[d]$}

Now consider corner $d_{ij}^w$ of $a_{ij}^w$, where $w \sim \sigma \in P_n^k$.

The matrix $d_{ij}^w$ has 1 in $(\alpha,\beta)$ row and $(\gamma,\delta)$ column
if and only if there exists $m \in \{1,\ldots,n\}$ such that
$S_\alpha S_\beta^* = S_i^* w S_\gamma S_m S_m^* S_\delta^* w^* S_j$,
if and only if there exists some $m$ such that
\begin{align}
(i,\alpha) & = \sigma(\gamma,m), \\
(j,\beta) & = \sigma(\delta,m) \nonumber.
\end{align}
Each row of $d_{ij}^w$ can have once 1 or be all 0's.
Summing over all $d_{ij}^w$, $i,j=1,\ldots,n$, each column has 1 in at most $n$ places
(possibly less).

\medskip

Let $\W_n^{k-1}$ be the union of the set of off-diagonal elements of $W_n^{k-1} \times W_n^{k-1}$
and $\{\dagger\}$, where $\dagger$ is a symbol not in $W_n^{k-1} \times W_n^{k-1}$.
Define mappings $f_{ij}^w: \W_n^{k-1} \to \W_n^{k-1}$ as
\begin{equation}
f_{ij}^w(\alpha,\beta) = (\gamma,\delta)
\end{equation}
if the entry of $d_{ij}^w$ in row $(\alpha,\beta)$ and column $(\gamma,\delta)$ is 1,
and as
\begin{equation}
f_{ij}^w(\alpha,\beta) = \dagger
\end{equation}
if the $(\alpha,\beta)$ row of $d_{ij}^w$ consists of all 0's.
In the latter case we think of $f_{ij}^w$ as ``annihilating'' $(\alpha,\beta)$.
Also, we put $f_{ij}(\dagger)=\dagger$ for all $i,j$.

Then $d_{ij}^w d_{rs}^w$ corresponds to $f_{rs}^w \circ f_{ij}^w$.
Again, in the sequel we drop the superscript $w$ when no confusion may arise.

We omit an easy proof of the following proposition.
\begin{lemma}
Let $w \in \P_n^k$. Then matrices $\{[d_{ij}^w]:i,j=1,\ldots,n\}$ generate a
nilpotent ring if and only if all sufficiently long composition products of mappings
$\{f_{ij}\;|\; i,j=1,\ldots,n\}$ have ranges consisting of the single element $\dagger$.
\end{lemma}
\begin{lemma}\label{dnil}
Let $w \in \P_n^k$. Then matrices $\{[d_{ij}^w]:i,j=1,\ldots,n\}$ generate a
nilpotent ring if and only there exists a partial order $\leq$ on $\W_n^{k-1}$ such that:
\begin{itemize}
\item[(i)]
The only minimal element with respect to $\leq$ is $\dagger$.
\item[(ii)]
For every $(\alpha,\beta) \in \W_n^{k-1}$ and all $i,j = 1, \ldots, n$,
\begin{equation}\label{fij}
f_{ij}(\alpha,\beta) \leq (\alpha,\beta).
\end{equation}
\end{itemize}
\end{lemma}
\begin{proof}
Suppose that the ring generated by $\{[d_{ij}^w]:i,j=1,\ldots,n\}$ is nilpotent.
Define a binary relation $\leq$ in $\W_n^{k-1}$ by (\ref{fij}) and take its
reflexive and transitive closure. Suppose for a moment that $(\alpha,\beta)\neq
(\gamma,\delta)$ but both $(\alpha,\beta)\leq(\gamma,\delta)$ and $(\gamma,\delta)
\leq(\alpha,\beta)$. Then, by definition of $\leq$, there are sequences $i_1,
\ldots,i_k$ and $j_1,\ldots j_k$ such that $f_{i_1 j_1}\circ \ldots \circ f_{i_k j_k}
(\alpha,\beta)=(\alpha,\beta)$. But then all composition powers of $f_{i_1 j_1}
\circ \ldots \circ f_{i_k j_k}$ have $(\alpha,\beta)$ in their range, a contradiction.

Conversely, suppose that there is a partial order $\leq$ on $\W_n^{k-1}$ satisfying
condition (ii) above. Then, by counting, each sufficiently long composition
product of mappings $\{f_{ij}\}$ has range consisting of a single element,
which is minimal for $\leq$. By (i), this element must be $\dagger$.
\end{proof}

Let $w\in\P_n^k$. We define inductively a nested sequence of subsets $\Psi_m^w$
of $\W_n^{k-1}$, as follows:
\begin{align}
\Psi_0^w & = \{\dagger\}, \\
\Psi_{m+1}^w & = \{(\alpha,\beta)\in\W_n^{k-1} \ | \
f_{ij}(\alpha,\beta) \in \Psi_m^w, \,i,j=1,\dots,n\} \cup \{\dagger\} \ .
\end{align}
We omit an easy proof of the following proposition.
\begin{proposition}\label{condd}
There exists a relation $\leq$ satisfying conditions of Lemma \ref{dnil}
if and only if
\begin{equation}\label{dcond}
\bigcup_m \Psi_m^w = \W_n^{k-1} \ .
\end{equation}
\end{proposition}

\subsection{A characterization of automorphisms in $\lambda(\P_n)^{-1}$}

From Theorem \ref{main}, Lemma \ref{nsfi} and Lemma \ref{dnil} we
obtain the following.
\begin{corollary}\label{charaut}
Let $w \in \P_n^k$. Then $\lambda_w \in {\rm Aut}(\O_n)$ if and only if the
following two conditions are satisfied:
\begin{enumerate}
\item There exists a partial order on $W_n^{k-1} \times W_n^{k-1}$
satisfying conditions of Lemma \ref{nsfi};
\item There exists a partial order on $\W_n^{k-1}$ satisfying conditions of Lemma \ref{dnil}.
\end{enumerate}
\end{corollary}

\section{Applications of labeled trees to automorphisms of $\O_2$}\label{otwo}

If $w\in\P_2^k$ then the labeled trees associated with $f_1^w$ and $f_2^w$
have the following properties:
\begin{itemize}
\item $\alpha$ receives two edges in $f_i^w$ if and only if $\alpha$ receives no
edges in $f_{3-i}^w$;
\item $\alpha$ receives one edge in $f_i^w$ if and only if $\alpha$ receives
one edge in $f_{3-i}^w$.
\end{itemize}
It follows that the numbers of leaves (0-receivers) on both trees are identical and
coincide with the number of 2-receivers (including the root) on these trees.
In such a case we say these two (unlabeled) trees are matched.

Given $w\in\P_2^k$ with corresponding functions $f_1^w$, $f_2^w$ and fixed $i\in\{1,2\}$,
we define
\begin{equation}\label{prestab}
G(f_i^w) := \{\sigma\in P_2^{k-1} \ | \ \sigma f_i^w\sigma^{-1}=f_i^w\},
\end{equation}
and call it the stabilizing group of $f_i^w$. Let $T$ be the unlabeled rooted tree
corresponding to $f_i^w$. If $\phi\in P_2^{k-1}$ then we have $G(f_i^w)
\cong G(\phi f_i^w \phi^{-1})$, through the map $\sigma\mapsto\phi\sigma\phi^{-1}$.
Thus the groups $G(f_i^w)$ do not depend on the choice of labels and we have
\begin{equation}\label{stab}
G(f_i^w)\cong{\rm Aut}(T),
\end{equation}
where ${\rm Aut}(T)$ is the automorphism group of the unlabeled rooted tree $T$.
Of course, a similar construction can be carried over for any $n$.

\subsection{Case of $\P_2^2$}

This case has been already well studied.
There are precisely four permutations in $\P_2^2$ yielding
automorphisms of $\O_2$. If $F := S_1 S_2^* + S_2 S_1^* \in \F_2^1$
denotes the flip-flop self-adjoint unitary, the four automorphisms are
${\rm id}, \lambda_F, {\rm Ad}(F) = \lambda_{\varphi(F)F}= \lambda_{F \varphi(F)},
{\rm Ad}(F)\lambda_F = \lambda_{\varphi(F)}$.
They form in ${\rm Aut}(\O_2)$ a copy of Klein's four-group.
In ${\rm Out}(\O_2)$, they give ${\mathbb Z}_2$ with
nontrivial generator the class of Archbold's flip-flop (Bogolubov)
automorphism $\lambda_F$, see e.g. \cite{Kaw1,Kaw2}.

Our labeled tree approach gives all these results with almost no effort at all.
The only pair of labeled trees satisfying Lemma \ref{nsfi} is
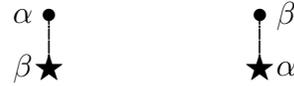
\[ \beginpicture
\setcoordinatesystem units <0.7cm,0.7cm>
\setplotarea x from 0 to 5, y from -1 to 0.1

\put {$\bullet$} at 0 0
\put {$\bigstar$} at 0 -1

\setlinear
\plot 0 0 0 -1 /
\put {$\alpha$} at -0.5 0
\put {$\beta$} at -0.5 -1

\put {$\bullet$} at 4 0
\put {$\bigstar$} at 4 -1

\setlinear
\plot 4 0 4 -1 /
\put {$\beta$} at 4.5 0
\put {$\alpha$} at 4.5 -1

\endpicture \]
Each is realized by 4 permutations and there are 2 such labelings.
Thus there are $2! \cdot 2^2 = 2 \cdot 4 = 8$ permutations in $P_2^2$ yielding
elements of ${\rm Aut}(\D_2)$. Of these 8 only 4 give automorphisms of $\O_2$.

\subsection{Case of $\P_2^3$}

Only two graphs are possible (each self-dual), namely

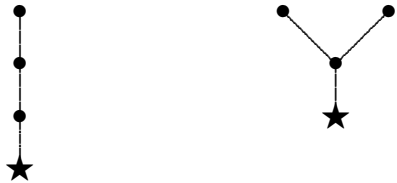
\[ \beginpicture
\setcoordinatesystem units <0.7cm,0.7cm>
\setplotarea x from 0 to 5, y from -1 to 1

\put {$\bullet$} at 0 1
\put {$\bullet$} at 0 0
\put {$\bullet$} at 0 -1
\put {$\bigstar$} at 0 -2

\setlinear
\plot 0 1 0 0 /
\plot 0 0 0 -1 /
\plot 0 -1 0 -2 /

\put {$\bullet$} at 5 1
\put {$\bullet$} at 7 1
\put {$\bullet$} at 6 0
\put {$\bigstar$} at 6 -1

\setlinear
\plot 5 1 6 0 /
\plot 7 1 6 0 /
\plot 6 0 6 -1 /

\endpicture \]

However, there is no labeling of the first graph which yields  correct
partial order $\leq$ on pairs. So only the second graph remains.
The only possible labeling satisfying conditions of Lemma \ref{nsfi} is
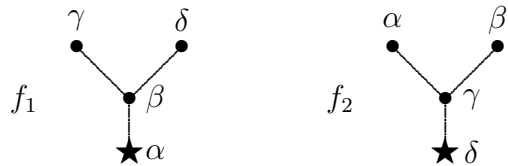
\[ \beginpicture
\setcoordinatesystem units <0.7cm,0.7cm>
\setplotarea x from 4 to 9, y from -1 to 1

\put {$\bullet$} at -1 1
\put {$\bullet$} at 1 1
\put {$\bullet$} at 0 0
\put {$\bigstar$} at 0 -1

\setlinear
\plot -1 1 0 0 /
\plot 1 1 0 0 /
\plot 0 0 0 -1 /
\put {$f_1$} at -2 0
\put {$\gamma$} at -1 1.5
\put {$\delta$} at 1 1.5
\put {$\beta$} at 0.5 0
\put {$\alpha$} at 0.5 -1

\put {$\bullet$} at 5 1
\put {$\bullet$} at 7 1
\put {$\bullet$} at 6 0
\put {$\bigstar$} at 6 -1

\setlinear
\plot 5 1 6 0 /
\plot 7 1 6 0 /
\plot 6 0 6 -1 /
\put {$f_2$} at 4 0
\put {$\alpha$} at 5 1.5
\put {$\beta$} at 7 1.5
\put {$\gamma$} at 6.5 0
\put {$\delta$} at 6.5 -1

\endpicture \]

Given a pair of labeled trees as above,
there are $2^4$ permutations $\sigma \in P_2^3$ yielding that pair.
There are $4!$ possible choices of labels. Hence, there are
\begin{equation}
4! \cdot 2^4 = 24 \cdot 16 = 324
\end{equation}
permutations in $P_2^3$ satisfying the conditions of Lemma \ref{nsfi}
and thus yielding elements of ${\rm Aut}(\D_2)$.

Then considering 16 permutations giving rise to a fixed labeling, as above,
one finds that only two of them satisfy the conditions of Lemma \ref{dnil}.
Thus, taking into account the action of inner automorphisms corresponding
to permutations in $P_2^2$, we see that there are exactly 48 automorphisms
of ${\mathcal O}_2$ corresponding to permutations in $P_2^3$. These are precisely
the ones inner equivalent to the identity or the flip-flop.
Thus, very surprisingly, among $8! = 40,320$ endomorphisms of ${\mathcal O}_2$
from $\lambda(\P_2^3)$ the only outer automorphism is the familiar flip-flop.
This is in stark contrast with the case of Cuntz algebras $\O_n$ with $n\geq3$,
where numerous new outer automorphisms appear already in $\lambda(\P_n^2)$
(see tables in Section 6.2, below).

Despite a large scale of the problem, our techniques allowed us to
obtain these results through easy and straightforward pen and paper
calculations. These were further confirmed through brute force computer
calculation based on the direct approach of Section 6.1, below.

\subsection{Case of $\P_2^4$}\label{ptwofour}

We begin by determining the number of automorphisms in $\lambda(\P_2^4)$.

\begin{theorem}\label{howmany}
We have
\begin{align*}
\# \{\lambda_w \ | \ w\in\P_2^4 \mbox{ and } \lambda_w|_{\D_2}
\in{\rm Aut}(\D_2)\} & = 8! \cdot 2^8 \cdot 17 = 175,472,640 \ , \\
\# \{\lambda_w \ | \ w\in\P_2^4 \mbox{ and } \lambda_w
\in{\rm Aut}(\O_2)\} & = 8! \cdot 14 = 564,480 \ .
\end{align*}
Thus in $\lambda(\P_2^4)^{-1}$ there are exactly $14$ representatives
of distinct inner equivalence classes.
\end{theorem}
\begin{proof}
There are exactly 23 directed rooted trees (unlabeled) with 8 vertices
satisfying our conditions (i.e. each vertex other than the root emits one edge
and receives maximum 2 edges, the root
is a minimal element and receives one edge from a different vertex). A computer
calculation shows that there are only 3 matched pairs of such trees admitting labelings
satisfying conditions of Proposition \ref{condb}. These are: $T_A-T_A$,
$T_A-T_J$ and $T_J-T_A$, where $T_A$ and $T_J$ are as follows (downward oriented):
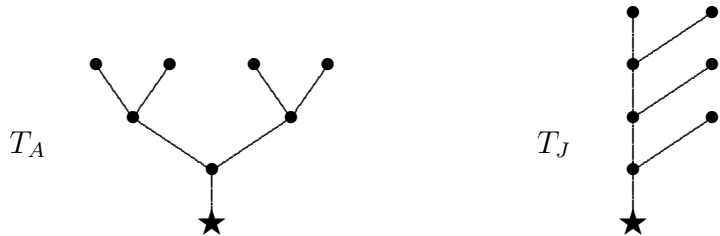
\[ \beginpicture
\setcoordinatesystem units <0.7cm,0.7cm>
\setplotarea x from 4 to 9, y from -1.5 to 3

\put {$\bullet$} at -2.2 2
\put {$\bullet$} at -0.8 2
\put {$\bullet$} at -3.8 2
\put {$\bullet$} at -5.2 2
\put {$\bullet$} at -4.5 1
\put {$\bullet$} at -1.5 1
\put {$\bullet$} at -3 0
\put {$\bigstar$} at -3 -1

\setlinear

\plot -0.8 2 -1.5 1 /
\plot -2.2 2 -1.5 1 /
\plot -5.2 2 -4.5 1 /
\plot -3.8 2 -4.5 1 /
\plot -4.5 1 -3 0 /
\plot -1.5 1 -3 0 /
\plot -3 0 -3 -1 /
\put {$T_A$} at -6.5 0.5

\put {$\bullet$} at 6.5 3
\put {$\bullet$} at 6.5 2
\put {$\bullet$} at 6.5 1
\put {$\bullet$} at 5 3
\put {$\bullet$} at 5 2
\put {$\bullet$} at 5 1
\put {$\bullet$} at 5 0
\put {$\bigstar$} at 5 -1

\setlinear

\plot 6.5 3 5 2 /
\plot 6.5 2 5 1 /
\plot 6.5 1 5 0 /
\plot 5 3 5 -1 /
\put {$T_J$} at 3.5 0.5

\endpicture \]
We fix arbitrarily labels on one of the trees in each pair, taking it
to be $T_J$ in the second and third case. Then computer calculation shows
the following numbers of labelings of the other tree which satisfy (\ref{bcond}):
40 for the pair $T_A-T_A$ and 12 for each of the other two pairs. The groups of
automorphisms of the rooted trees $T_A$ and $T_J$ have 8 and 2 elements,
respectively. Thus, taking into account that each pair of labeled trees under
consideration is realized by $2^8$ distinct permutations, and factoring in the action
of $8!$ inner automorphisms (which permute the labels simultaneously on both trees),
we obtain the following number of distinct permutations in $P_2^4$ giving rise
to automorphisms of the diagonal:
$$ 2^8\cdot\frac{8!}{|{\rm Aut}(T_A)|}\cdot 40 + 2\cdot 2^8 \cdot
\frac{8!}{|{\rm Aut}(T_J)|}\cdot 12 = 2^8\cdot 8!\cdot 17 = 175,472,640. $$
Then a computer calculation shows that among these permutations there are only
$8!\cdot 14=564,480$ satisfying (\ref{dcond})
and thus yielding automorphisms of $\O_2$. Dividing out $8!$ inner
automorphisms from level $3$, we finally get $14$ inner equivalence classes
of automorphisms in $\lambda(\P_2^4)^{-1}$.
\end{proof}

Our next goal is to describe explicitly representatives of inner equivalence
classes from $\lambda(\P_2^4)^{-1}$ and to find some infinite subgroups of
${\rm Out}(\O_2)$ generated by them.

We begin by considering two permutations $A$ and $B$ of the set $W_2^4$ given
respectively by

\begin{equation*}
\begin{array}{cccc}
A(1211)=1211 &  A(1212)=1212 &  A(1222)=1222 &  A(1221)=1221 \\
A(1121)=1121 &  A(1122)=1122 &  A(1111)=1112 &  A(1112)=1111 \\
A(2222)=2111 &  A(2221)=2121 &  A(2211)=2112  & A(2212)=2122 \\
A(2122)=2222  & A(2121)=2221 &  A(2112)=2212  &  A(2111)=2211
\end{array}
\end{equation*}

\begin{equation*}
\begin{array}{cccc}
B(1211)=1211 &  B(1212)=1212 &  B(1222)=1222 &  B(1221)=1221 \\
B(1121)=1121 &  B(1122)=1122 &  B(1111)=1112 &  B(1112)=1111 \\
B(2122)=2111 &  B(2121)=2112 &  B(2211)=2121 &  B(2212)=2122 \\
B(2222)=2212 &  B(2221)=2221 &  B(2112)=2222 &  B(2111)=2211
\end{array}
\end{equation*}

Note that the first two rows of these two permutations are identical.
That is, $A(1***)=B(1***)$.
And of the first eight arguments, six are fixed points.
The labeled trees corresponding to $A$ are:

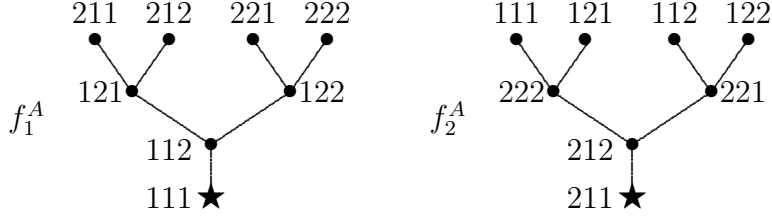
\[ \beginpicture
\setcoordinatesystem units <0.7cm,0.7cm>
\setplotarea x from 4 to 9, y from -1.5 to 1

\put {$\bullet$} at -2.2 2
\put {$\bullet$} at -0.8 2
\put {$\bullet$} at -3.8 2
\put {$\bullet$} at -5.2 2
\put {$\bullet$} at -4.5 1
\put {$\bullet$} at -1.5 1
\put {$\bullet$} at -3 0
\put {$\bigstar$} at -3 -1

\setlinear

\plot -0.8 2 -1.5 1 /
\plot -2.2 2 -1.5 1 /
\plot -5.2 2 -4.5 1 /
\plot -3.8 2 -4.5 1 /
\plot -4.5 1 -3 0 /
\plot -1.5 1 -3 0 /
\plot -3 0 -3 -1 /
\put {$f_1^A$} at -6.5 0.5
\put {$222$} at -0.8 2.5
\put {$221$} at -2.2 2.5
\put {$211$} at -5.2 2.5
\put {$212$} at -3.8 2.5
\put {$121$} at -5.1 1
\put {$122$} at -0.9 1
\put {$112$} at -3.8 -0.1
\put {$111$} at -3.8 -1

\put {$\bullet$} at 5.8 2
\put {$\bullet$} at 7.2 2
\put {$\bullet$} at 4.1 2
\put {$\bullet$} at 2.8 2
\put {$\bullet$} at 3.5 1
\put {$\bullet$} at 6.5 1
\put {$\bullet$} at 5 0
\put {$\bigstar$} at 5 -1

\setlinear

\plot 7.2 2 6.5 1 /
\plot 5.8 2 6.5 1 /
\plot 2.8 2 3.5 1 /
\plot 4.2 2 3.5 1 /
\plot 3.5 1 5 0 /
\plot 6.5 1 5 0 /
\plot 5 0 5 -1 /
\put {$f_2^A$} at 1.5 0.5
\put {$122$} at 7.2 2.5
\put {$112$} at 5.8 2.5
\put {$111$} at 2.8 2.5
\put {$121$} at 4.2 2.5
\put {$222$} at 2.9 1
\put {$221$} at 7.1 1
\put {$212$} at 4.2 -0.1
\put {$211$} at 4.2 -1

\endpicture \]

In the sequel, for notational convenience, we equip $W_2^k$
with the reversed lexicographic order
and enumerate its elements as $\{1, 2, \ldots, 2^k\}$ accordingly.
Then, the permutations $A$ and $B$
above correspond to
$A = (1,9)(2,4,10,12,14,16)(6,8)$ and
$B = (1,9)(2,4,6,10,16,12,14)$.
With a slight abuse of notation we also denote simply by
$A$ and $B$ the associated unitaries and by
$\lambda_A$ and $\lambda_B$
the corresponding endomorphisms of $\O_2$.

Using Corollary \ref{charaut} one can verify that $\lambda_A$ and $\lambda_B$
are automorphisms of $\O_2$. In fact, these permutations were found
through pen and paper calculation based on Corollary \ref{charaut}.
One checks by computer calculation based on Section 6.1 that the inverses of the
automorphisms $\lambda_A$ and $\lambda_B$ are induced by unitaries in $\P_2^7$.
\begin{proposition}
In ${\rm Out}(\O_2)$, one has
$$\lambda_F \lambda_A \lambda_F = \lambda_A^{-1} = \lambda_B \ . $$
\end{proposition}
\begin{proof}
One has ${\rm Ad}(z) \lambda_A \lambda_B = {\rm id}$, where $z \in \P_2^6$ is given by
\begin{align*}
z \sim & (2, 4, 8)(3, 7, 15)(5, 13, 29)(9, 25)(10, 12) \\
    & (18, 20, 24)(19, 23)(26, 28)(34, 36, 40) \\
    & (35, 39, 47)(37, 45)(42, 44)(50, 52, 56)(51, 55)(58, 60).
\end{align*}
Also, one has ${\rm Ad}(y) \lambda_F \lambda_A = \lambda_B \lambda_F$, where
$y \sim ( 1, 3, 5, 7) (2, 4, 8) \in P_2^3$.
\end{proof}

For reader's convenience, in Appendix \ref{applA} we provide the action
of $\lambda_A$ on diagonal projections $P_\alpha$'s
with $|\alpha| \leq 5$.

\begin{lemma}\label{lap}
With the above notation, for each word $\tilde{\mu}$ there exist words
$\nu_1,\nu_2$ with $|\nu_i| = |\tilde\mu| + 1$ such that
\begin{align*}
\lambda_A(P_{\tilde\mu 211}) & = P_{\nu_1 211} + P_{\nu_2 222} \ , \\
\lambda_A(P_{\tilde\mu 212}) & = P_{\nu_1 212} + P_{\nu_2 221} \ .
\end{align*}
Furthermore, if $\alpha$ is a word which ends neither with 211 nor with 212
then there is a word $\beta$ such that $|\alpha|=|\beta|$ and
$\lambda_A(P_\alpha)=P_\beta$.
\end{lemma}
\begin{proof}
We proof the first claim by induction on $|\tilde\mu|$.
If $|\tilde\mu| \leq 2$ these relations are verified by direct computation.
Now let us suppose that $\tilde\mu = (\mu_1,\ldots,\mu_l)$ and $l \geq 3$. Then
\begin{align*}
\lambda_A(P_{\tilde\mu 211}) & = \lambda_A(P_{\mu_1 \ldots \mu_l 211})
= A^* S_{\mu_1} \lambda_A(P_{\mu_2 \ldots \mu_l 211}) S_{\mu_1}^* A \\
& = A^* S_{\mu_1} (P_{\tilde\nu_1 211} + P_{\tilde\nu_2 222}) S_{\mu_1}^* A \\
& = A^*(P_{\mu_1 \tilde\nu_1 211} + P_{\mu_1 \tilde\nu_2 222})A \\
& = P_{\cdots 211} + P_{\cdots 222}
\end{align*}
where in the second line we have used the induction hypothesis
and in the last line we have used the fact that $|\mu_1 \tilde\nu_1| = |\mu_1 \tilde\nu_2| \geq 4$.
The other relation can be handled similarly.

The proof of the second claim proceeds by induction on $|\alpha|$. For $|\alpha| \leq 3$
this follows from the table in Appendix \ref{applA}. For the inductive step we notice
that there exist two unitaries $u_1,u_2$ in ${\mathcal F}_2^3$ such that
$\lambda_A(S_i)=S_i u_i$, $i=1,2$. Thus, we have $\lambda_A(P_{i\alpha})=
\lambda_A(S_i)\lambda_A(P_\alpha)\lambda_A(S_i)^* = S_i u_i P_\beta u_i^* S_i^* =
P_{i\mu}$ for some word $\mu$ with $|\mu|=|\alpha|$.
\end{proof}

\begin{proposition}\label{aorder}
$\lambda_A$ has infinite order
in ${\rm Out}(\O_2)$.
\end{proposition}

 \begin{proof}
It is a consequence of Lemma \ref{lap} that $\lambda_A$ has infinite order in ${\rm Aut}(\O_2)$.
To see this, fix some $\tilde\mu$.
If some power of $\lambda_A$ were the identity then, using the relations in Lemma \ref{lap},
one should have that $P_{\tilde\mu 211}$ is a sum of subprojections including one of the form $P_{\rho 211}$.
But then $P_{\rho 212}$ should also be a subprojection of $P_{\tilde\mu 211}$.
On the other hand,
by the same relations $P_{\rho 212}$ should be subprojection of $P_{\tilde\mu 212}$ and
thus orthogonal to $P_{\tilde\mu 211}$,
contradiction.

Now it follows from implication (1) $\Rightarrow$ (2) of \cite[Theorem 6]{Sz} that
$\lambda_A$ has infinite order in ${\rm Out}(\O_2)$.
\end{proof}

\begin{corollary}
The subgroup of ${\rm Out}(\O_2)$ generated by $\lambda_A$ and $\lambda_F$ is
isomorphic to the infinite dihedral group ${\mathbb Z} \rtimes {\mathbb Z}_2$.
\end{corollary}

Let $J$ be a transposition in $P_2^4$ which exchanges $2112$ with $2212$ (and
fixes all other elements of $W_2^4$):
$$ J(2112)=2212 \;\;\; \mbox{and} \;\;\; J(2212)=2112. $$
The labeled trees corresponding to $J$ are:
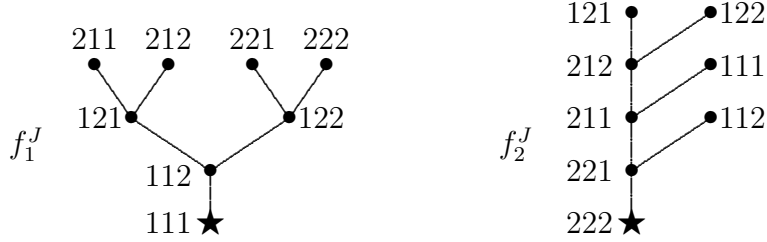
\[ \beginpicture
\setcoordinatesystem units <0.7cm,0.7cm>
\setplotarea x from 4 to 9, y from -1.5 to 3

\put {$\bullet$} at -2.2 2
\put {$\bullet$} at -0.8 2
\put {$\bullet$} at -3.8 2
\put {$\bullet$} at -5.2 2
\put {$\bullet$} at -4.5 1
\put {$\bullet$} at -1.5 1
\put {$\bullet$} at -3 0
\put {$\bigstar$} at -3 -1

\put {$222$} at -0.8 2.5
\put {$221$} at -2.2 2.5
\put {$211$} at -5.2 2.5
\put {$212$} at -3.8 2.5
\put {$121$} at -5.1 1
\put {$122$} at -0.9 1
\put {$112$} at -3.8 -0.1
\put {$111$} at -3.8 -1

\setlinear

\plot -0.8 2 -1.5 1 /
\plot -2.2 2 -1.5 1 /
\plot -5.2 2 -4.5 1 /
\plot -3.8 2 -4.5 1 /
\plot -4.5 1 -3 0 /
\plot -1.5 1 -3 0 /
\plot -3 0 -3 -1 /
\put {$f_1^J$} at -6.5 0.5

\put {$\bullet$} at 6.5 3
\put {$\bullet$} at 6.5 2
\put {$\bullet$} at 6.5 1
\put {$\bullet$} at 5 3
\put {$\bullet$} at 5 2
\put {$\bullet$} at 5 1
\put {$\bullet$} at 5 0
\put {$\bigstar$} at 5 -1

\put {$121$} at 4.2 3
\put {$212$} at 4.2 2
\put {$211$} at 4.2 1
\put {$221$} at 4.2 0
\put {$222$} at 4.2 -1
\put {$122$} at 7.1 3
\put {$111$} at 7.1 2
\put {$112$} at 7.1 1

\setlinear

\plot 6.5 3 5 2 /
\plot 6.5 2 5 1 /
\plot 6.5 1 5 0 /
\plot 5 3 5 -1 /
\put {$f_2^J$} at 2.8 0.5

\endpicture \]
With a slight abuse of notation, we denote by $J$ the associated unitary
and by $\lambda_J$ the corresponding endomorphism of $\O_2$. One checks that
\begin{equation}
\lambda_J^2={\rm id}.
\end{equation}
Clearly (see Example \ref{bogolubov}), the two trees corresponding to the identity
in $P_2^4$ are both of type $T_A$. Likewise, both trees corresponding to the
flip-flop $\lambda_F$ are also of type $T_A$. Since $f_2^J$ is of type $T_J\neq T_A$,
it follows from Proposition \ref{outclasses} that $\lambda_J$ is an outer automorphism
of $\O_2$ not inner equivalent to the flip-flop. Incidentally, outerness of
$\lambda_J$ can also be derived from \cite{MaTo}, since $\lambda_J(S_1)=S_1$.
\begin{proposition}
Automorphisms $\lambda_F$ and $\lambda_J$ generate a subgroup of ${\rm Out}(\O_2)$
isomorphic to the free product ${\mathbb Z}_2*{\mathbb Z}_2$.
\end{proposition}
\begin{proof}
The proof is very similar to the argument of Proposition \ref{aorder} and Lemma
\ref{lap}, so we only sketch the main idea.

At first one shows by induction on word length that for each word $\mu$ there exists
a word $\nu$ with $|\nu|=|\mu|+2$ such that
\begin{align*}
\lambda_F\lambda_J(P_{\mu22}) & = P_{\nu22}+\sum_i P_{\gamma_i} \ ,\\
\lambda_F\lambda_J(P_{\mu21}) & = P_{\nu21}+\sum_j P_{\zeta_j} \ ,
\end{align*}
with $\gamma_i,\zeta_j$ words of lengths not greater than $|\nu|+2$. This implies
that automorphism $\lambda_F\lambda_J$ has infinite order. Consequently, it
has an infinite order in ${\rm Out}(\O_2)$, and the claim follows.
\end{proof}

Let $G$ be a 3-cycle in $W_2^4$ such that
$$ G(1112)=1122, \;\;\; G(1122)=1222, \;\; \mbox{ and } \;\; G(1222)=1112. $$
That is, in the shorthand notation, $G=(9,13,15)$. The trees corresponding to $G$ are:
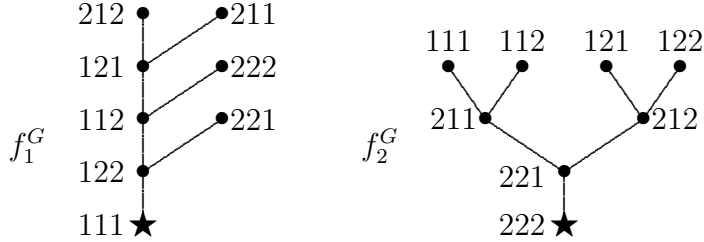
\[ \beginpicture
\setcoordinatesystem units <0.7cm,0.7cm>
\setplotarea x from 4 to 9, y from -1.5 to 3

\put {$\bullet$} at 5.8 2
\put {$\bullet$} at 7.2 2
\put {$\bullet$} at 4.2 2
\put {$\bullet$} at 2.8 2
\put {$\bullet$} at 3.5 1
\put {$\bullet$} at 6.5 1
\put {$\bullet$} at 5 0
\put {$\bigstar$} at 5 -1

\put {$122$} at 7.2 2.5
\put {$121$} at 5.8 2.5
\put {$111$} at 2.8 2.5
\put {$112$} at 4.2 2.5
\put {$211$} at 2.9 1
\put {$212$} at 7.1 1
\put {$221$} at 4.2 -0.1
\put {$222$} at 4.2 -1

\setlinear

\plot 7.2 2 6.5 1 /
\plot 5.8 2 6.5 1 /
\plot 2.8 2 3.5 1 /
\plot 4.2 2 3.5 1 /
\plot 3.5 1 5 0 /
\plot 6.5 1 5 0 /
\plot 5 0 5 -1 /

\put {$f_2^G$} at 1.5 0.5

\put {$\bullet$} at -1.5 3
\put {$\bullet$} at -1.5 2
\put {$\bullet$} at -1.5 1
\put {$\bullet$} at -3 3
\put {$\bullet$} at -3 2
\put {$\bullet$} at -3 1
\put {$\bullet$} at -3 0
\put {$\bigstar$} at -3 -1

\put {$212$} at -3.8 3
\put {$121$} at -3.8 2
\put {$112$} at -3.8 1
\put {$122$} at -3.8 0
\put {$111$} at -3.8 -1
\put {$211$} at -0.9 3
\put {$222$} at -0.9 2
\put {$221$} at -0.9 1

\setlinear

\plot -1.5 3 -3 2 /
\plot -1.5 2 -3 1 /
\plot -1.5 1 -3 0 /
\plot -3 3 -3 -1 /

\put {$f_1^G$} at -5.2 0.5

\endpicture \]
One checks that
\begin{equation}
\lambda_G^6={\rm id}
\end{equation}
but none of $\lambda_G$, $\lambda_G^2$, $\lambda_G^3$ is inner. Also note that
$\lambda_G(S_2)=S_2$.

Taking into account the results of this subsection and
considering the convolution multiplication and Lemma \ref{l:2.3}
(and preferably helped by a computer), one verifies
the following theorem.

\begin{theorem}
The following automorphisms give a complete list of representatives of distinct
classes in ${\rm Out}(\O_2)$ appearing in $\lambda(\P_2^4)^{-1}$:
\begin{align*}
& \{ {\rm id}, \; \lambda_F \}, \\
& \{ \lambda_A, \; \lambda_A \lambda_F, \; \lambda_F \lambda_A, \; \lambda_F \lambda_A \lambda_F \}, \\
& \{ \lambda_J, \; \lambda_J \lambda_F, \; \lambda_F \lambda_J, \; \lambda_F \lambda_J \lambda_F \}, \\
& \{ \lambda_G, \; \lambda_G \lambda_F, \; \lambda_F \lambda_G, \; \lambda_F \lambda_G \lambda_F \}.
\end{align*}
\end{theorem}

\section{Computations and tables}\label{computables}

\subsection{Inverse pairs of localized automorphisms}\label{pme}

In this short section we gather together a few facts about pairs of
unitaries in some finite matrix algebras giving rise to
automorphisms of $\O_n$ that are inverses of each other.
We also briefly discuss interesting algebraic equations such
unitaries must satisfy. These equations provide
useful background for the considerations in Section 3
(e.g. Theorem \ref{main}, Corollary \ref{bound}).
They have also been useful for concrete computations, e.g. in
computing explicitly the inverse of $\lambda_A$
in Section 4, filling the tables of the following subsection,
and in the search of square-free automorphisms.

Hereafter,
for any unitary $u \in \O_n$ and a positive integer $k$, we set
\begin{equation}
u_k := u^* \varphi(u^*) \cdots \varphi^{k-1}(u^*).
\end{equation}
Notice that $u_k \varphi^k (u_h) = u_{k+h}$. In this subsection, symbol
$u_k^*$ should always be understood as $(u_k)^*$.

So let us suppose that $U \in \F_n^k$, $V \in \F_n^h$ are unitaries such that
$$\lambda_U \lambda_V = {\rm id} = \lambda_V \lambda_U  \ , $$
i.e. $U\lambda_U(V) = 1 = V\lambda_V(U)$.\footnote{
Since $\lambda_U$ and $\lambda_V$ are injective, one identity implies the other.
Also, up to replacing $k$ and $h$ with
$k \vee h$
there would be no loss of generality in assuming that $k = h$, however as
the inverse of an automorphism induced by a unitary in a
matrix algebra might very well be induced by a unitary in a
larger matrix algebra it seems convenient to allow this more flexible asymmetric formulation.
It is worth stressing that, given $k$, the subset of unitaries $U$'s in $\F_n^k$
such that $\lambda_U^{-1}$ (exists and) is still induced by a unitary in $\F_n^k$
is definitely smaller than the set of unitaries such that $\lambda_U^{-1}$ is induced
by a unitary in some $\F_n^h$.
An a priori bound for $h$ as a function of $n,k$ is provided by Corollary \ref{bound}.}
Then we readily obtain the coupled system of matrix equations
\begin{equation}\label{coupled}
U_h V U_h^* = U^*, \quad V_k U V_k^* = V^*,
\end{equation}
where both $U_h$ and $V_k$ are in $\F_n^{h + k-1}$.
In passing, observe that the second equation is independent of the level $h$
for which $V\in\F_n^h$.

In practical situations, one is faced with the converse problem.
Starting with some $U \in \F_n^k$, one might not know the precise value of $h$,
let alone if the corresponding $V$ exists at all.
It turns out that solutions (for $V$) of equations (\ref{coupled}) imply
invertibility of $\lambda_U$. The following proposition combined with Corollary \ref{bound}
gives an algorithmic procedure for finding these solutions. We omit an elementary proof.
\begin{proposition}
Let $U$ be a unitary in $\F_n^k$ and suppose that $U_h^* U^* U_h \in \F_n^h$ for some $h$.
Then $\lambda_U$ is invertible and $\lambda_U^{-1}=\lambda_V$ with $V := U_h^* U^* U_h$.
\end{proposition}

In particular, given a unitary $U \in \F_n^k$, one has $\lambda_U^2 = {\rm id}$
(i.e., $U = V$)
if and only if $U\lambda_U(U) = 1$, if and only if $U_k U U_k^* = U^*$.

\medskip
Finally, we present yet another computational strategy for determining invertibility of
endomorphism $\lambda_U$ and finding its inverse. Again, we omit an elementary proof
of the following proposition.
\begin{proposition}\label{equ}
Let $U$ and $V$ be unitaries in $\F_n^k$ and $\F_n^h$, respectively, satisfying
equations (\ref{coupled}). Then $U$ is a solution of the following polynomial matrix equation
\begin{equation}\label{necU}
(U_r^* U^* U_r)_r U (U_r^* U^* U_r)_r^* = U_r^* U U_r \ ,
\end{equation}
where $r$ can be taken as maximum of $k$ and $h$.

Conversely, given $r$, every solution $U \in \F_n^r$ of equation (\ref{necU})
gives rise to an automorphism $\lambda_U$ of $\O_n$,
with inverse induced by $V:=U_r^* U^* U_r$.
\end{proposition}
\begin{remark}
\rm The strategy of Proposition \ref{equ} is to find all pairs satisfying (\ref{coupled})
by solving equations of the form (\ref{necU}) for all values of $r$.
Implicitly, by solving such an equation, we predict $V$ to take a
particular form, namely $V=U_r^* U^* U_r$. However, we do not assume $V\in\F_n^r$.
In fact, $V$ automatically belongs to $\F_n^{2r-1}$. Combining this with equations (\ref{coupled})
we obtain an additional relation $U$ must satisfy, namely
$U_r^* U^* U_r = U_{2r-1}^* U^* U_{2r-1}$.
\end{remark}

We find it rather intriguing that
in the case of permutation unitaries
the polynomial matrix equations (\ref{necU}) turn out
to be equivalent to the tree related conditions of Corollary \ref{charaut}.

\subsection{Tabulated results}

In this section, we collect our results about automorphisms $\lambda(\P_n)^{-1}$
of the Cuntz algebras in the form of tables.
They provide solutions to several enumeration problems.

In the first table, we provide the number $N_n^k$ of all such
automorphisms of $\O_n$ at level $k$ (i.e. in $\lambda(\P_n^k)^{-1}$),
for small values of $n$ and $k$.
In the second table, we plot the number $C_n^k$ of classes modulo inner ones.
Of course, we have
$$N_n^k = n^{k-1}! \; C_n^k \ . $$
The last table contains numbers $sf_n^k$ of square-free
automorphisms in $\lambda(\P_n^k)^{-1}$.

\bigskip

\begin{center}
$N_n^k$: \quad
\begin{tabular}{|l||l|l|l|}
\hline
$k \setminus n$ & 2 & 3 & 4\\
\hline\hline
1 & 2  &  6 & 24 \\
\hline
2 & 4  &  576 & 5,771,520 \\
\hline
3 & 48 &  & \\
\hline
4 & 564,480 &  & \\
\hline
\end{tabular}
\end{center}

\begin{center}
$C_n^k$: \quad
\begin{tabular}{|l||l|l|l|}
\hline
$k \setminus n$ & 2 & 3 & 4 \\
\hline\hline
1 & 2  & 6 & 24 \\
\hline
2 & 2 & 96  & 240,480  \\
\hline
3 & 2 & & \\
\hline
4 & 14  &  & \\
\hline
\end{tabular}
\end{center}

\begin{center}
$sf_n^k$: \quad
\begin{tabular}{|l||l|l|l|}
\hline
$k \setminus n$ & 2 & 3 & 4 \\
\hline\hline
1 & 2 & 4 & 10 \\
\hline
2 & 4 & 52 & 2,032 \\
\hline
3 & 20 & & \\
\hline
4 & 1,548 &  & \\
\hline
\end{tabular}
\end{center}

\medskip

These figures have been obtained through combination of all the techniques developed in
this article and large scale computer calculations. In particular, both labeled tree approach of
Corollary \ref{charaut} and algebraic equation approach of Section 6.1 have been
used. To give the reader an idea of the scale of the problem and difficulties involved
let us just mention that computation of $N_4^2$ (and thus $C_4^2$) took about
70 processor days.

\section{Concluding remarks}

If $n\geq3$ then the image of $\lambda(\P_n^2)^{-1}$ in ${\rm Out}(\O_n)$ contains
${\mathbb Z}_3*{\mathbb Z}_2$ (see \cite{Sz}) and thus it is non-amenable. In
the case of $\O_2$ we still do not known if the group $\lambda(\P_n)^{-1}$ (and
its image in ${\rm Out}(\O_2)$) is amenable or not. It would be interesting to find
the lowest level $k$ (if any) for which $\lambda(\P_2^k)^{-1}$ is non-amenable.
Our results show that $k$ must be at least $4$, and this question can perhaps
be settled by determining the group generated by $\lambda(\P_2^4)^{-1}$.

Going beyond automorphisms of $\O_n$ preserving the $UHF$-subalgebra, one may
pose the question if any aspects of the theory developed in the present article
can be extended to $\lambda({\mathcal S}_n)^{-1}$. This is certainly far from obvious and
undoubtedly a very challenging task. Even the inner part of $\lambda({\mathcal S}_n)^{-1}$,
that is the group ${\mathcal S}_n$ itself, is non-amenable and has a very complicated structure.
In fact, in the case of $\O_2$, it
contains a copy of ${\mathbb Z}_3*{\mathbb Z}_2$ whose action on the diagonal results
in the crossed product isomorphic to $\O_2$ \cite{Sp}. Finding a criterion of invertibility of
endomorphisms from $\lambda({\mathcal S}_n)$ should certainly be regarded as an
important first step.

Our labeled tree approach allows for relatively easy construction of certain special
automorphisms of the diagonal ${\mathcal D}_n$, which in turn give rise to dynamical systems
on the Cantor set. They certainly deserve further investigations. In particular,
a question arises if they may result in minimal dynamical systems. Also, their relation
with other better known classes of symbolic dynamical systems is worth elucidation.

Of course, the polynomial matrix equations of Subsection \ref{pme} apply to
arbitrary unitaries in the algebraic part of ${\mathcal F}_n$ and not only to
permutation matrices. Therefore, they can be used for finding other families of
automorphisms of $\O_n$ with localized inverses. It is to be expected that new
interesting classes of automorphisms different from the much studied quasi-free ones
will be found this way. It is also intriguing to investigate the algebraic variety
defined by these equations. At present, we are not aware of occurences
of these equations outside the realm of Cuntz algebras but we would not be surprised
if such instances were found.

Going beyond the Cuntz algebras, it is natural to expect that parts of our analysis
may be extended to more general graph $C^*$-algebras. Every graph algebra admits a
gauge action of the circle group whose fixed point algebra is $AF$ (e.g. see \cite{Ra}).
For a large class of graph algebras, a Cartan subalgebra with totally disconnected
spectrum is contained in this core $AF$-subalgebra. Then, for such graph algebras,
one should be able to say much about automorphisms preserving both the core $AF$ and
the Cartan subalgebra along the lines of \cite{Cun2}, \cite{Sz} and the present article.

\section{Appendix}\label{applA}

\begin{tabular}{|l||l|l|l|} \hline
$P_\alpha$ & $\lambda_A(P_\alpha)$ & $\lambda_G(P_\alpha)$ & $\lambda_J(P_\alpha)$ \\ \hline\hline
$P_1$ & $P_1$ & $P_1$ & $P_1$ \\
$P_2$ & $P_2$ & $P_2$ & $P_2$ \\ \hline
$P_{11}$ & $P_{11}$ & $P_{111}+P_{1121}+P_{1222}$ & $P_{11}$ \\
$P_{12}$ & $P_{12}$ & $P_{1122}+P_{121}+P_{1221}$ & $P_{12}$ \\
$P_{21}$ & $P_{22}$ & $P_{21}$ & $P_{2111}+P_{212}+P_{2212}$ \\
$P_{22}$ & $P_{21}$ & $P_{22}$ & $P_{2112}+P_{2211}+P_{222}$ \\ \hline
$P_{111}$ & $P_{111}$ & $P_{1111}+P_{12221}+P_{11122}$ & $P_{111}$ \\
$P_{112}$ & $P_{112}$ & $P_{11121}+P_{1121}+P_{12222}$ & $P_{112}$ \\
$P_{121}$ & $P_{122}$ & $P_{121}$ & $P_{12111}+P_{1212}+P_{12212}$ \\
$P_{122}$ & $P_{212}$ & $P_{1122}+P_{1221}$ & $P_{12112}+P_{12211}+P_{1222}$ \\
$P_{211}$ & $P_{2211}+P_{2222}$ & $P_{2111}+P_{21121}+P_{21222}$ & $P_{2111}+P_{2212}$ \\
$P_{212}$ & $P_{2212}+P_{2221}$ & $P_{21122}+P_{2121}+P_{21221}$ & $P_{212}$ \\
$P_{221}$ & $P_{212}$ & $P_{221}$ & $P_{2112}+P_{22111}+P_{22212}$ \\
$P_{222}$ & $P_{211}$ & $P_{222}$ & $P_{22112}+P_{22211}+P_{2222}$ \\ \hline
$P_{1111}$ & $P_{1112}$ & $P_{11111}+P_{111122}+P_{111221}$ & $P_{1111}$ \\
$P_{1112}$ & $P_{1111}$ & $P_{111121}+P_{111222}+P_{12221}$ & $P_{1112}$ \\
$P_{1121}$ & $P_{1122}$ & $P_{1121}$ & $P_{112111}+P_{11212}+P_{112212}$ \\
$P_{1122}$ & $P_{1121}$ & $P_{12222}+P_{11121}$ & $P_{112112}+P_{112211}+P_{11222}$ \\
$P_{1211}$ & $P_{12211}+P_{12222}$ & $P_{12111}+P_{121121}+P_{121222}$ & $P_{12111}+P_{12212}$ \\
$P_{1212}$ & $P_{12212}+P_{12221}$ & $P_{121122}+P_{12121}+P_{121221}$ & $P_{1212}$ \\
$P_{1221}$ & $P_{1212}$ & $P_{1221}$ & $P_{12112}+P_{122111}+P_{122212}$ \\
$P_{1222}$ & $P_{1211}$ & $P_{1122}$ & $P_{122112}+P_{122211}+P_{12222}$ \\
$P_{2111}$ & $P_{2222}$ & $P_{21111}+P_{211122}+P_{212221}$ & $P_{2111}$ \\
$P_{2112}$ & $P_{2211}$ & $P_{211121}+P_{21121}+P_{212222}$ & $P_{2212}$ \\
$P_{2121}$ & $P_{2212}$ & $P_{2121}$ & $P_{212111}+P_{21212}+P_{212212}$ \\
$P_{2122}$ & $P_{2221}$ & $P_{21122}+P_{21221}$ & $P_{212112}+P_{212211}+P_{21222}$ \\
$P_{2211}$ & $P_{21211}+P_{21222}$ & $P_{22111}+P_{221121}+P_{221222}$ & $P_{22111}+P_{22212}$ \\
$P_{2212}$ & $P_{21212}+P_{21221}$ & $P_{221122}+P_{22121}+P_{221221}$ & $P_{2112}$ \\
$P_{2221}$ & $P_{2112}$ & $P_{2221}$ & $P_{22112}+P_{222111}+P_{222212}$ \\
$P_{2222}$ & $P_{2111}$ & $P_{2222}$ & $P_{222112}+P_{222211}+P_{22222}$ \\ \hline
\end{tabular}

\vspace{3mm}\noindent
\begin{tabular}{|l||l|l|l|} \hline
$P_\alpha$ & $\lambda_A(P_\alpha)$ & $\lambda_G(P_\alpha)$ & $\lambda_J(P_\alpha)$ \\ \hline\hline
$P_{11111}$ &  $P_{11122}$ & $P_{111111}+P_{1111122}+P_{1111221}$ & $P_{11111}$ \\
$P_{11112}$ &  $P_{11121}$ & $P_{1111121}+P_{1111222}+P_{111221}$ & $P_{11112}$ \\
$P_{11121}$ &  $P_{11112}$ & $P_{12221}$ & $P_{1112111}+P_{111212}+P_{1112212}$ \\
$P_{11122}$ &  $P_{11111}$ & $P_{111121}+P_{111222}$ & $P_{1112112}+P_{1112211}+P_{111222}$ \\
$P_{11211}$ &  $P_{112211}+P_{112222}$ & $P_{112111}+P_{1121121}+P_{1121222}$ & $P_{112111}+P_{112212}$ \\
$P_{11212}$ &  $P_{112212}+P_{112221}$ & $P_{1121122}+P_{112121}+P_{1121221}$ & $P_{11212}$ \\
$P_{11221}$ &  $P_{11212}$ & $P_{11121}$ & $P_{112112}+P_{1122111}+P_{1122212}$ \\
$P_{11222}$ &  $P_{11211}$ & $P_{12222}$ & $P_{1122112}+P_{1122211}+P_{112222}$ \\
$P_{12111}$ &  $P_{12222}$ & $P_{121111}+P_{1211122}+P_{1212221}$ & $P_{12111}$ \\
$P_{12112}$ &  $P_{12211}$ & $P_{1211121}+P_{121121}+P_{1212222}$ & $P_{12212}$ \\
$P_{12121}$ &  $P_{12212}$ & $P_{12121}$ & $P_{1212111}+P_{121212}+P_{1212212}$ \\
$P_{12122}$ &  $P_{12221}$ & $P_{121122}+P_{121221}$ & $P_{1212112}+P_{1212211}+P_{121222}$ \\
$P_{12211}$ &  $P_{121211}+P_{121222}$ & $P_{122111}+P_{1221121}+P_{1221222}$ & $P_{122111}+P_{122212}$ \\
$P_{12212}$ &  $P_{121212}+P_{121221}$ & $P_{1221122}+P_{122121}+P_{1221221}$ & $P_{12112}$ \\
$P_{12221}$ &  $P_{12112}$ & $P_{11221}$ & $P_{122112}+P_{1222111}+P_{1222212}$ \\
$P_{12222}$ &  $P_{12111}$ & $P_{11222}$ & $P_{1222112}+P_{1222211}+P_{122222}$ \\
$P_{21111}$ &  $P_{22222}$ & $P_{211111}+P_{2111122}+P_{2111221}$ & $P_{21111}$ \\
$P_{21112}$ &  $P_{22221}$ & $P_{2111121}+P_{2111222}+P_{212221}$ & $P_{21112}$ \\
$P_{21121}$ &  $P_{22112}$ & $P_{21121}$ & $P_{2212111}+P_{221212}+P_{2212212}$ \\
$P_{21122}$ &  $P_{22111}$ & $P_{212222}+P_{211121}$ & $P_{2212112}+P_{2212211}+P_{221222}$ \\
$P_{21211}$ &  $P_{221211}+P_{221222}$ & $P_{212111}+P_{2121121}+P_{2121222}$ & $P_{212111}+P_{212212}$ \\
$P_{21212}$ &  $P_{221212}+P_{221221}$ & $P_{2121122}+P_{212121}+P_{2121221}$ & $P_{21212}$ \\
$P_{21221}$ &  $P_{22212}$ & $P_{21221}$ & $P_{212112}+P_{2122111}+P_{2122212}$ \\
$P_{21222}$ &  $P_{22211}$ & $P_{21122}$ & $P_{2122112}+P_{2122211}+P_{212222}$ \\
$P_{22111}$ &  $P_{21222}$ & $P_{221111}+P_{2211122}+P_{2212221}$ & $P_{22111}$ \\
$P_{22112}$ &  $P_{21211}$ & $P_{2211121}+P_{221121}+P_{2212222}$ & $P_{22212}$ \\
$P_{22121}$ &  $P_{21212}$ & $P_{22121}$ & $P_{2112111}+P_{211212}+P_{2112212}$ \\
$P_{22122}$ &  $P_{21221}$ & $P_{221122}+P_{221221}$ & $P_{2112112}+P_{2112211}+P_{211222}$ \\
$P_{22211}$ &  $P_{211211}+P_{211222}$ & $P_{222111}+P_{2221121}+P_{2221222}$ & $P_{222111}+P_{222212}$ \\
$P_{22212}$ &  $P_{211212}+P_{211221}$ & $P_{2221122}+P_{222121}+P_{2221221}$ & $P_{22112}$ \\
$P_{22221}$ &  $P_{21112}$ & $P_{22221}$ & $P_{222112}+P_{2222111}+P_{2222212}$ \\
$P_{22222}$ &  $P_{21111}$ & $P_{22222}$ & $P_{2222112}+P_{2222211}+P_{222222}$ \\ \hline
\end{tabular}

\newpage \noindent
Roberto Conti\\
Mathematics, School of Mathematical and Physical Sciences \\
University of Newcastle, Callaghan, NSW 2308, Australia
\\ E-mail: Roberto.Conti@newcastle.edu.au \\

\smallskip \noindent
Wojciech Szyma{\'n}ski\\
Mathematics, School of Mathematical and Physical Sciences \\
University of Newcastle, Callaghan, NSW 2308, Australia
\\ E-mail: Wojciech.Szymanski@newcastle.edu.au \\

\end{document}